\documentclass[fleqn]{article}

\usepackage{amsmath}
\usepackage{amssymb}
\usepackage{amsthm}

\usepackage[dvips]{graphicx} 

\title{
  The $2$-bridge knots of up to $16$ crossings
}

\author{David De Wit%
  \footnote{%
    Department of Mathematics, The University of Queensland,
    4072, Brisbane, Australia.
    \newline
    \texttt{Dr\_David\_De\_Wit@yahoo.com.au}
    }
}

\begin{document}

\maketitle

\begin{abstract}
  \noindent
  For any given number of crossings $c$, there exists a formula to determine
  the number of $2$-bridge knots of $c$ crossings, and indeed it is a simple
  matter to actually construct presentations of these knots. However, the
  determination of whether a given (prime) knot is a $2$-bridge knot remains a
  nontrivial exercise, and we have no procedure to determine bridge numbers
  more generally. Herein, we identify the $2$-bridge knots within the
  Hoste--Thistlethwaite--Weeks tables of prime knots of up to $16$ crossings
  by an exhaustive search of a larger set of $2$-bridge knots. As the unknot
  is the only knot with bridge number $1$, this yields a lower bound of $3$
  for the bridge numbers of the remaining knots.
\end{abstract}


\section{Introduction}

Recall that the $2$-bridge knots are a class of alternating prime knots. Most,
but not all, are chiral, and as oriented knots, all are invertible. Although
it is straightforward to determine a catalogue of $2$-bridge knots, the
determination of whether a given (prime) knot is a $2$-bridge knot, is
nontrivial. In principle, one takes the $2$-fold branched cover of the knot,
then computes the characteristic (geometric) decomposition (by the orbifold
theorem, we know such a decomposition exists). If the $2$-fold cover is a lens
space, then the knot is $2$-bridge, otherwise not. Herein, we do not pursue
this algorithm, rather, we identify the $2$-bridge knots in the
Hoste--Thistlethwaite--Weeks tables by exhaustively examining a larger
catalogue of $2$-bridge knots.

Adopting the notation of Ernst and Sumners~\cite{ErnstSumners:1987}, denote by
$TK_c$ (respectively $TK^*_c$) the number of $2$-bridge knots of $c$ crossings
when a knot and its reflection are regarded as nondistinct (respectively
distinct). Also let $ATK_c$ denote the number of achiral $2$-bridge knots of
$c$ crossings, so we have: $TK^*_c = 2 TK_c - ATK_c$. As $2$-bridge knots can
be described purely by sequences of integers (see below), they may be counted
combinatorially. Indeed,~\cite{ErnstSumners:1987} presents the following
formulae for $c\geqslant 3$, which in Table~\ref{table:TKc} we evaluate for
$c\leqslant 16$. Note that there are no achiral alternating knots, $2$-bridge
or otherwise, with odd crossing numbers~\cite{HosteThistlethwaiteWeeks:1998}.

\begin{eqnarray*}
  TK^*_c
  & = &
  \left\{
    \begin{array}{ll}
      {\scriptstyle \frac{1}{3}} (2^{c-2}-1)             & c ~\mathrm{even} \\
      {\scriptstyle \frac{1}{3}} (2^{c-2}+2^{(c-1)/2})   & c \equiv 1~(\mathrm{mod}~4) \\
      {\scriptstyle \frac{1}{3}} (2^{c-2}+2^{(c-1)/2}+2) & c \equiv 3~(\mathrm{mod}~4)
    \end{array}
  \right.
  \\
  \\
  ATK_c
  & = &
  \left\{
    \begin{array}{ll}
      {\scriptstyle \frac{1}{3}} (2^{(c-2)/2}+1) \qquad \quad & c \equiv 0~(\mathrm{mod}~4) \\
      {\scriptstyle \frac{1}{3}} (2^{(c-2)/2}-1)              & c \equiv 2~(\mathrm{mod}~4) \\
      0                                                       & c ~\mathrm{odd}.
    \end{array}
  \right.
\end{eqnarray*}

\renewcommand{\arraystretch}{1.5}
\renewcommand{\tabcolsep}{4pt}

\begin{table}[htbp]
  \tiny
  \begin{centering}
  \begin{tabular}{r|*{14}{r}|r}
   \multicolumn{1}{c}{} & \multicolumn{14}{c}{$c$} \\
                  & 3 & 4 & 5 & 6 & 7 & 8 & 9 & 10 & 11 & 12 & 13 & 14 & 15 & 16 &total \\
   \hline
     $TK^*_c$ & 2 & 1 & 4 & 5 & 14 & 21 & 48 & 85 & 182 & 341 & 704 & 1365 & 2774 & 5461 & 11007 \\
      $ATK_c$ & . & 1 & . & 1 & . & 3 & . & 5 & . & 11 & . & 21 & . & 43 & 85 \\
       $TK_c$ & 1 & 1 & 2 & 3 & 7 & 12 & 24 & 45 & 91 & 176 & 352 & 693 & 1387 & 2752 & 5546 \\
      $K^A_c$ & 1 & 1 & 2 & 3 & 7 & 18 & 41 & 123 & 367 & 1288 & 4878 & 19536 & 85263 & 379799 & 491327
  \end{tabular}
  \caption{%
    Numbers of $2$-bridge and alternating prime knots.
  }%
  \label{table:TKc}
  \end{centering}
\end{table}

The notation is compatible with the following conventions. Let $K_c$ denote
the number of prime knots of $c$ crossings, and then denote by $K^A_c$
(respectively $K^N_c$) the number of alternating (respectively nonalternating)
prime knots of $c$ crossings. Next, let $c^A_i$ (respectively $c^N_i$) denote
the $i$th alternating (respectively nonalternating) prime knot (type) of $c$
crossings in the HTW tables of prime knots of up to $16$
crossings~\cite{HosteThistlethwaiteWeeks:1998}. Note that the convention is to
use a superscript ``$A$'' to denote an alternating knot, whilst the ``$A$'' in
$ATK_c$ denotes ``achiral''. The presence of the superscript decorating the
crossing number identifies the knot as from the HTW tables, and facilitates
the concurrent use of the classical undecorated Alexander--Briggs notation;
for instance $5^A_1$ and $5_2$ denote the same knot type.
Dowker--Thistlethwaite codes and many other data associated with the knots in
the HTW tables may be accessed via the program \textsc{Knotscape} (version
1.01), as mentioned in~\cite{HosteThistlethwaiteWeeks:1998}.

Now, each $2$-bridge knot may be expressed via a presentation known as
\emph{Conway's normal form}. By restricting the degrees of freedom used in
choosing this presentation, each $2$-bridge knot may be uniquely expressed,
indeed this remains true when a knot and its reflection are regarded as
distinct. Importantly however, the canonical presentation is not generally
minimal, that is, the presentation will usually have more crossings than the
prime knot type to which it corresponds. So, by exhaustively enumerating and
identifying sufficiently many canonical presentations, we may identify all the
$2$-bridge knots within the HTW tables. This approach to enumerating the
$2$-bridge knots is less mathematically elegant than that
of~\cite{ErnstSumners:1987}, but is easy to implement.

Herein, we compile a list of the $2$-bridge knots of up to $16$ crossings by
first enumerating all such canonical Conway presentations of up to $28$
crossings. To each presentation we assign an orientation, and then use this to
deduce an appropriate DT code. This process is automated using
\textsc{Mathematica}. We then use the \texttt{Locate} (that is, the
``\texttt{Locate in Table}'') function of \textsc{Knotscape} to identify which
elements of our list of DT codes correspond to which (alternating) prime knots
within the HTW tables.

Now let $P$ be a Conway presentation identified as being a knot of type $T$
described by a DT code within the tables. In general, $T$ will be chiral, and
it is of interest to ask to which chirality of $T$ our presentation $P$
corresponds. However, any DT code necessarily corresponds to both a knot and
its reflection (the book by Adams~\cite[p38]{Adams:1994} shows how). That is,
the HTW tables do not provide us with canonical representatives of knot types
in the way that classical pictorial knot tables do, so we may not identify $P$
with a particular chiral knot, only a knot \emph{type}. In contrast, the DT
codes of the HTW tables do prescribe orientations whereas classical tables do
not. However, $2$-bridge knots are invertible, so we may ignore this issue.
All that said, \textsc{Knotscape} happily evaluates various polynomial
invariants of knot types described by DT codes; so its routines make an opaque
choice which determines a chirality.

Apart from its intrinsic interest, we have compiled our list of $2$-bridge
knots to facilitate our enquiries into the properties of the Links--Gould
invariant $LG^{2,1}$~\cite{DeWit:2000,DeWitKauffmanLinks:1999a}. Currently,
skein relations sufficient to determine the value of this invariant for any
arbitrary link remain unknown, and the only general method for its evaluation
is a computationally-expensive state model method. This method is currently
computationally feasible for links for which braid presentations of at most
$5$ strings are available; to date we have been able to evaluate $LG^{2,1}$
this way only for $37,547$ of the $1,701,936$ knots in the HTW tables.
However, recent work by Ishii~\cite{Ishii:2004b} has determined a formula for
$LG^{2,1}$ for $2$-bridge links described by canonical Conway presentations.
So, by identifying $2$-bridge knots within the HTW tables, we may evaluate
$LG^{2,1}$ for many prime knots for which the state model method remains
infeasible. Braids for the HTW knots may be obtained via the use of the
program \textsc{K2K} by Imafuji and Ochiai~\cite{ImafujiOchiai:2002}; this
program is a \textsc{Mathematica} interface to the \textsc{C} program
\textsc{KnotTheoryByComputer}. We observe that the string indices of
\textsc{K2K}-generated braids corresponding to our $2$-bridge knots are
generally greater than $5$; indeed this reflects their inherent complexity.
The results of the application of Ishii's formula to evaluations of $LG^{2,1}$
for $2$-bridge knots from the HTW tables, together with further material on
the current state of evaluations of $LG^{2,1}$, appear
in~\cite{DeWitLinks:WhereLG21Fails}.

Our intensive use of \textsc{Knotscape} in this manner demonstrates the
robustness of its algorithms for reducing Dowker--Thistlethwaite codes. This
is significant as it is known that these algorithms are necessarily incomplete.

We mention that~\cite{ErnstSumners:1987} also includes formulae for the
numbers $TL^*_c$, $TL_c$ and $ATL_c$ of $2$-component $2$-bridge links. We
have not pursued the identification of these links here as the existing tables
of multicomponent prime links are not as extensive as those of the true knots.
More generally, the methodology applied to identify the $2$-bridge knots can be
applied to identify the torus and pretzel knots. Although we currently have no
formulae for $LG^{2,1}$ for knots of these classes, relatively efficient
versions of state model algorithms can be implemented for them. However, these
knots, particularly the pretzels, are comparatively rare in the HTW tables,
and we reach a point of diminishing returns.


\section{Conway presentations for $2$-bridge links}

\begin{figure}[ht]
  \begin{centering}
  \begin{picture}(340,120)
 \qbezier(10,70)(1,71)(0,80)
 \qbezier(0,80)(1,89)(10,90)
 \qbezier(10,110)(1,110.5)(0,115)
 \qbezier(0,115)(1,119.5)(10,120)
 \qbezier(10,120)(90,120)(170,120)

 \qbezier(10,110)(15,110)(17.5,100)
 \qbezier(17.5,100)(20,90)(25,90)
 \qbezier(25,110)(20,110)(18.5,104)
 \qbezier(16.5,96)(15,90)(10,90)
 \multiput(27,100)(3,0){3}{\circle*{1}}
 \qbezier(35,110)(40,110)(42.5,100)
 \qbezier(42.5,100)(45,90)(50,90)
 \qbezier(50,110)(45,110)(43.5,104)
 \qbezier(41.5,96)(40,90)(35,90)
 \put(30,85){\makebox(0,0){$a_1$}}
 \qbezier(10,70)(30,70)(50,70)

 \qbezier(50,70)(55,70)(57.5,80)
 \qbezier(57.5,80)(60,90)(65,90)
 \qbezier(65,70)(60,70)(58.5,76)
 \qbezier(56.5,84)(55,90)(50,90)
 \multiput(67,80)(3,0){3}{\circle*{1}}
 \qbezier(75,70)(80,70)(82.5,80)
 \qbezier(82.5,80)(85,90)(90,90)
 \qbezier(90,70)(85,70)(83.5,76)
 \qbezier(81.5,84)(80,90)(75,90)
 \put(70,95){\makebox(0,0){$a_2$}}
 \qbezier(50,110)(70,110)(90,110)

 \qbezier(90,110)(95,110)(97.5,100)
 \qbezier(97.5,100)(100,90)(105,90)
 \qbezier(105,110)(100,110)(98.5,104)
 \qbezier(96.5,96)(95,90)(90,90)
 \multiput(107,100)(3,0){3}{\circle*{1}}
 \qbezier(115,110)(120,110)(122.5,100)
 \qbezier(122.5,100)(125,90)(130,90)
 \qbezier(130,110)(125,110)(123.5,104)
 \qbezier(121.5,96)(120,90)(115,90)
 \put(110,85){\makebox(0,0){$a_3$}}
 \qbezier(90,70)(110,70)(130,70)

 \qbezier(130,70)(135,70)(137.5,80)
 \qbezier(137.5,80)(140,90)(145,90)
 \qbezier(145,70)(140,70)(138.5,76)
 \qbezier(136.5,84)(135,90)(130,90)
 \multiput(147,80)(3,0){3}{\circle*{1}}
 \qbezier(155,70)(160,70)(162.5,80)
 \qbezier(162.5,80)(165,90)(170,90)
 \qbezier(170,70)(165,70)(163.5,76)
 \qbezier(161.5,84)(160,90)(155,90)
 \put(150,95){\makebox(0,0){$a_4$}}
 \qbezier(130,110)(150,110)(170,110)

 \multiput(173,70)(3,0){3}{\circle*{1}}
 \multiput(173,90)(3,0){3}{\circle*{1}}
 \multiput(173,110)(3,0){3}{\circle*{1}}
 \multiput(173,120)(3,0){3}{\circle*{1}}
 \multiput(196,70)(3,0){3}{\circle*{1}}
 \multiput(196,90)(3,0){3}{\circle*{1}}
 \multiput(196,110)(3,0){3}{\circle*{1}}
 \multiput(196,120)(3,0){3}{\circle*{1}}

 \qbezier(205,110)(210,110)(212.5,100)
 \qbezier(212.5,100)(215,90)(220,90)
 \qbezier(220,110)(215,110)(213.5,104)
 \qbezier(211.5,96)(210,90)(205,90)
 \multiput(222,100)(3,0){3}{\circle*{1}}
 \qbezier(230,110)(235,110)(237.5,100)
 \qbezier(237.5,100)(240,90)(245,90)
 \qbezier(245,110)(240,110)(238.5,104)
 \qbezier(236.5,96)(235,90)(230,90)
 \put(225,85){\makebox(0,0){$a_n$}}
 \qbezier(205,70)(225,70)(245,70)

 \qbezier(245,90)(254,89)(255,80)
 \qbezier(255,80)(254,71)(245,70)
 \qbezier(245,120)(254,119.5)(255,115)
 \qbezier(255,115)(254,110.5)(245,110)
 \qbezier(205,120)(225,120)(245,120)
 \put(325,95){\makebox(0,0){$n$: odd}}

 \qbezier(10,0)(1,1)(0,10)
 \qbezier(0,10)(1,19)(10,20)
 \qbezier(10,40)(1,40.5)(0,45)
 \qbezier(0,45)(1,49.5)(10,50)
 \qbezier(10,50)(90,50)(170,50)

 \qbezier(10,40)(15,40)(17.5,30)
 \qbezier(17.5,30)(20,20)(25,20)
 \qbezier(25,40)(20,40)(18.5,34)
 \qbezier(16.5,26)(15,20)(10,20)
 \multiput(27,30)(3,0){3}{\circle*{1}}
 \qbezier(35,40)(40,40)(42.5,30)
 \qbezier(42.5,30)(45,20)(50,20)
 \qbezier(50,40)(45,40)(43.5,34)
 \qbezier(41.5,26)(40,20)(35,20)
 \put(30,15){\makebox(0,0){$a_1$}}
 \qbezier(10,0)(30,0)(50,0)

 \qbezier(50,0)(55,0)(57.5,10)
 \qbezier(57.5,10)(60,20)(65,20)
 \qbezier(65,0)(60,0)(58.5,6)
 \qbezier(56.5,14)(55,20)(50,20)
 \multiput(67,10)(3,0){3}{\circle*{1}}
 \qbezier(75,0)(80,0)(82.5,10)
 \qbezier(82.5,10)(85,20)(90,20)
 \qbezier(90,0)(85,0)(83.5,6)
 \qbezier(81.5,14)(80,20)(75,20)
 \put(70,25){\makebox(0,0){$a_2$}}
 \qbezier(50,40)(70,40)(90,40)

 \qbezier(90,40)(95,40)(97.5,30)
 \qbezier(97.5,30)(100,20)(105,20)
 \qbezier(105,40)(100,40)(98.5,34)
 \qbezier(96.5,26)(95,20)(90,20)
 \multiput(107,30)(3,0){3}{\circle*{1}}
 \qbezier(115,40)(120,40)(122.5,30)
 \qbezier(122.5,30)(125,20)(130,20)
 \qbezier(130,40)(125,40)(123.5,34)
 \qbezier(121.5,26)(120,20)(115,20)
 \put(110,15){\makebox(0,0){$a_3$}}
 \qbezier(90,0)(110,0)(130,0)

 \qbezier(130,0)(135,0)(137.5,10)
 \qbezier(137.5,10)(140,20)(145,20)
 \qbezier(145,0)(140,0)(138.5,6)
 \qbezier(136.5,14)(135,20)(130,20)
 \multiput(147,10)(3,0){3}{\circle*{1}}
 \qbezier(155,0)(160,0)(162.5,10)
 \qbezier(162.5,10)(165,20)(170,20)
 \qbezier(170,0)(165,0)(163.5,6)
 \qbezier(161.5,14)(160,20)(155,20)
 \put(150,25){\makebox(0,0){$a_4$}}
 \qbezier(130,40)(150,40)(170,40)

 \multiput(173,0)(3,0){3}{\circle*{1}}
 \multiput(173,20)(3,0){3}{\circle*{1}}
 \multiput(173,40)(3,0){3}{\circle*{1}}
 \multiput(173,50)(3,0){3}{\circle*{1}}
 \multiput(196,0)(3,0){3}{\circle*{1}}
 \multiput(196,20)(3,0){3}{\circle*{1}}
 \multiput(196,40)(3,0){3}{\circle*{1}}
 \multiput(196,50)(3,0){3}{\circle*{1}}

 \qbezier(205,40)(210,40)(212.5,30)
 \qbezier(212.5,30)(215,20)(220,20)
 \qbezier(220,40)(215,40)(213.5,34)
 \qbezier(211.5,26)(210,20)(205,20)
 \multiput(222,30)(3,0){3}{\circle*{1}}
 \qbezier(230,40)(235,40)(237.5,30)
 \qbezier(237.5,30)(240,20)(245,20)
 \qbezier(245,40)(240,40)(238.5,34)
 \qbezier(236.5,26)(235,20)(230,20)
 \put(225,15){\makebox(0,0){$a_{n-1}$}}
 \qbezier(205,0)(225,0)(245,0)

 \qbezier(245,0)(250,0)(252.5,10)
 \qbezier(252.5,10)(255,20)(260,20)
 \qbezier(260,0)(255,0)(253.5,6)
 \qbezier(251.5,14)(250,20)(245,20)
 \multiput(262,10)(3,0){3}{\circle*{1}}
 \qbezier(270,0)(275,0)(277.5,10)
 \qbezier(277.5,10)(280,20)(285,20)
 \qbezier(285,0)(280,0)(278.5,6)
 \qbezier(276.5,14)(275,20)(270,20)
 \put(265,25){\makebox(0,0){$a_n$}}
 \qbezier(245,40)(265,40)(285,40)

 \qbezier(285,40)(294,39)(295,30)
 \qbezier(295,30)(294,21)(285,20)
 \qbezier(305,20)(305,30)(305,40)
 \qbezier(285,50)(303,49)(305,40)
 \qbezier(305,20)(303,2)(285,0)
 \qbezier(205,50)(245,50)(285,50)
 \put(325,25){\makebox(0,0){$n$: even}}
\end{picture}
  \caption{%
    The $2$-bridge link presentation $C(a)\equiv C(a_1,\dots,a_n)$,
    illustrated for positive coefficients $a_i$. Reflect the crossings for
    negative $a_i$.
  }%
  \label{figure:ConwayPresentationForUnorientedTwoBridgeLinks}
  \end{centering}
\end{figure}
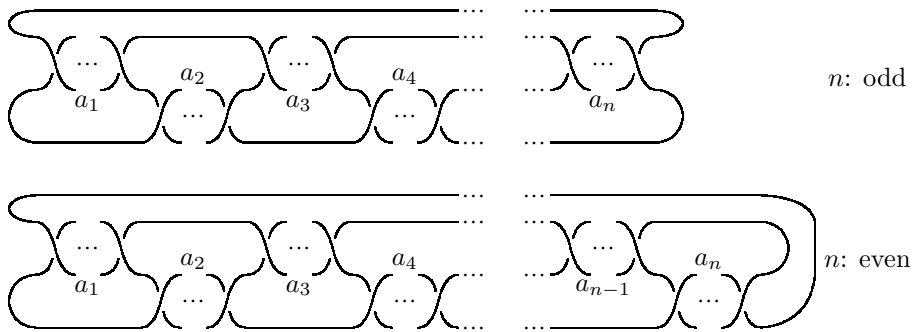

Abstracting material from the excellent books by
Kawauchi~\cite[pp21--26]{Kawauchi:1996} and
Murasugi~\cite[pp171--196]{Murasugi:1996}, we describe the $2$-bridge links
using the plait presentation known as \emph{Conway's normal form}, as depicted
in Figure~\ref{figure:ConwayPresentationForUnorientedTwoBridgeLinks}. In this
form, any $2$-bridge link may be described in terms of a presentation $C(a)$,
where $a$ is a finite sequence $a_1,\dots,a_n$ of nonzero integers $a_i$. In
the present context, we need not define the meaning of ``$2$-bridge''; we are
only interested in the $2$-bridge links as that class of links for which there
exist presentations of the normal form. Note the convention relating the signs
of the $a_i$ to those of the crossings, and observe that such a presentation
has either $1$ or $2$ components.

Now choose integers $\alpha$ and $\beta$ such that $\alpha>0$ and
$\gcd(\alpha,\beta)=1$ in the following continued fraction, referred to as the
\emph{slope} of the presentation $C(a)$:
\begin{eqnarray*}
  \frac{\alpha}{\beta}
  \triangleq
  a_1 + \frac{1}{a_2 + \frac{1}{a_3 + \cdots \frac{1}{a_n}}}.
\end{eqnarray*}
The presentation $C(a)$ is then said to be of $2$-bridge link type
$S(\alpha,\beta)$. After Murasugi~\cite[p180]{Murasugi:1996}, we denote the
continued fraction by $[a_1,\dots,a_n]$.

The $2$-bridge \emph{knots} are in fact precisely the \emph{rational} knots.
Conway presentations are equal (that is, they correspond to the same link)
when their slopes are equal; for example $C(2,-3)=C(1,1,2)=S(5,3)$, as
$[2,-3]=[1,1,2]=\frac{5}{3}$. More generally, $2$-component $2$-bridge links
$S(\alpha,\beta)$ and $S(\alpha',\beta')$ are equal if and only if
$\alpha=\alpha'$ and $\beta^{\pm 1} \equiv \beta'~(\mathrm{mod}~2\alpha)$, and
$2$-bridge knots $S(\alpha,\beta)$ and $S(\alpha',\beta')$ are equal if and
only if $\alpha=\alpha'$ and $\beta^{\pm 1} \equiv
\beta'~(\mathrm{mod}~\alpha)$. Next, observe that
$S(\alpha,\beta)^*=S(\alpha,-\beta)$; in fact, $S(\alpha,\beta)$ is achiral if
and only if $\beta^2=-1~(\mathrm{mod}~\alpha)$. Also observe that
$C(a)^*=C(-a)$, where by $-a$ we intend the sequence $-a_1,\dots,-a_n$.

Now, if for some $S(\alpha,\beta)$, we seek a representative presentation
$C(a)$, we may choose all the $a_i$ to be of the same sign as $\beta$. Such a
$C(a)$ is thus an alternating presentation, so $S(\alpha,\beta)$ is
necessarily an \emph{alternating} link. Moreover, any $S(\alpha,\beta)$ always
has a presentation $C(a)$ with only \emph{even} components $a_i$ (although the
$a_i$ are \emph{not} necessarily all of the same sign). For example,
$C(2,-3)=S(5,3)=S(5,2)=C(2,2)$. Knowing that any $2$-bridge link may be
expressed via such a presentation $C(a)$ with all $a_i$ even, we observe that
it is a knot if $n$ is even and a $2$-component link if $n$ is odd.

\begin{figure}[ht]
  \begin{centering}
  \begin{picture}(340,60)
 \qbezier(10,0)(1,1)(0,10)
 \qbezier(0,10)(1,19)(10,20)
 \qbezier(10,40)(1,40.5)(0,45)
 \qbezier(0,45)(1,49.5)(10,50)
 \qbezier(10,50)(90,50)(185,50)
 \put(28,50){\vector(-1,0){0}}

 \qbezier(10,40)(15,40)(17.5,30)
 \qbezier(17.5,30)(20,20)(25,20)
 \qbezier(25,40)(20,40)(18.5,34)
 \qbezier(16.5,26)(15,20)(10,20)
 \multiput(27,30)(3,0){3}{\circle*{1}}
 \qbezier(35,40)(40,40)(42.5,30)
 \qbezier(42.5,30)(45,20)(50,20)
 \qbezier(50,40)(45,40)(43.5,34)
 \qbezier(41.5,26)(40,20)(35,20)
 \put(30,15){\makebox(0,0){$a_1$}}
 \qbezier(10,0)(30,0)(55,0)
 \put(32,0){\vector(1,0){0}}
 \qbezier(50,20)(50,20)(55,20)
 \put(50,20){\vector(-1,0){0}}

 \qbezier(55,0)(60,0)(62.5,10)
 \qbezier(62.5,10)(65,20)(70,20)
 \qbezier(70,0)(65,0)(63.5,6)
 \qbezier(61.5,14)(60,20)(55,20)
 \multiput(72,10)(3,0){3}{\circle*{1}}
 \qbezier(80,0)(85,0)(87.5,10)
 \qbezier(87.5,10)(90,20)(95,20)
 \qbezier(95,0)(90,0)(88.5,6)
 \qbezier(86.5,14)(85,20)(80,20)
 \put(75,25){\makebox(0,0){$a_2$}}
 \qbezier(50,40)(75,40)(100,40)
 \put(77,40){\vector(1,0){0}}
 \qbezier(95,20)(95,20)(100,20)
 \put(95,20){\vector(-1,0){0}}

 \qbezier(100,40)(105,40)(107.5,30)
 \qbezier(107.5,30)(110,20)(115,20)
 \qbezier(115,40)(110,40)(108.5,34)
 \qbezier(106.5,26)(105,20)(100,20)
 \multiput(117,30)(3,0){3}{\circle*{1}}
 \qbezier(125,40)(130,40)(132.5,30)
 \qbezier(132.5,30)(135,20)(140,20)
 \qbezier(140,40)(135,40)(133.5,34)
 \qbezier(131.5,26)(130,20)(125,20)
 \put(120,15){\makebox(0,0){$a_3$}}
 \qbezier(95,0)(120,0)(145,0)
 \put(122,0){\vector(1,0){0}}
 \qbezier(140,20)(140,20)(145,20)
 \put(140,20){\vector(-1,0){0}}

 \qbezier(145,0)(150,0)(152.5,10)
 \qbezier(152.5,10)(155,20)(160,20)
 \qbezier(160,0)(155,0)(153.5,6)
 \qbezier(151.5,14)(150,20)(145,20)
 \multiput(162,10)(3,0){3}{\circle*{1}}
 \qbezier(170,0)(175,0)(177.5,10)
 \qbezier(177.5,10)(180,20)(185,20)
 \qbezier(185,0)(180,0)(178.5,6)
 \qbezier(176.5,14)(175,20)(170,20)
 \put(165,25){\makebox(0,0){$a_4$}}
 \qbezier(140,40)(165,40)(185,40)
 \put(167,40){\vector(1,0){0}}

 \multiput(188,0)(3,0){3}{\circle*{1}}
 \multiput(188,20)(3,0){3}{\circle*{1}}
 \multiput(188,40)(3,0){3}{\circle*{1}}
 \multiput(188,50)(3,0){3}{\circle*{1}}
 \multiput(211,0)(3,0){3}{\circle*{1}}
 \multiput(211,20)(3,0){3}{\circle*{1}}
 \multiput(211,40)(3,0){3}{\circle*{1}}
 \multiput(211,50)(3,0){3}{\circle*{1}}

 \qbezier(220,40)(225,40)(227.5,30)
 \qbezier(227.5,30)(230,20)(235,20)
 \qbezier(235,40)(230,40)(228.5,34)
 \qbezier(226.5,26)(225,20)(220,20)
 \multiput(237,30)(3,0){3}{\circle*{1}}
 \qbezier(245,40)(250,40)(252.5,30)
 \qbezier(252.5,30)(255,20)(260,20)
 \qbezier(260,40)(255,40)(253.5,34)
 \qbezier(251.5,26)(250,20)(245,20)
 \put(240,15){\makebox(0,0){$a_{n-1}$}}
 \qbezier(220,0)(240,0)(265,0)
 \put(243,0){\vector(1,0){0}}
 \qbezier(260,20)(260,20)(265,20)
 \put(260,20){\vector(-1,0){0}}

 \qbezier(265,0)(270,0)(272.5,10)
 \qbezier(272.5,10)(275,20)(280,20)
 \qbezier(280,0)(275,0)(273.5,6)
 \qbezier(271.5,14)(270,20)(265,20)
 \multiput(282,10)(3,0){3}{\circle*{1}}
 \qbezier(290,0)(295,0)(297.5,10)
 \qbezier(297.5,10)(300,20)(305,20)
 \qbezier(305,0)(300,0)(298.5,6)
 \qbezier(296.5,14)(295,20)(290,20)
 \put(285,25){\makebox(0,0){$a_n$}}
 \qbezier(260,40)(285,40)(305,40)
 \put(288,40){\vector(1,0){0}}

 \qbezier(305,40)(314,39)(315,30)
 \qbezier(315,30)(314,21)(305,20)
 \qbezier(325,20)(325,30)(325,40)
 \qbezier(305,50)(323,49)(325,40)
 \qbezier(325,20)(323,2)(305,0)
 \qbezier(220,50)(265,50)(305,50)
 \put(283,50){\vector(-1,0){0}}

\end{picture}
  \caption{%
    The oriented $2$-bridge knot even presentation $C(a)\equiv
    C(a_1,\dots,a_n)$, where both $n$ and each $a_i$ are nonzero even integers. The
    illustration depicts positive coefficients $a_i$; reflect the
    crossings for negative $a_i$.
  }%
  \label{figure:EvenConwayPresentationForOrientedTwoBridgeKnots}
  \end{centering}
\end{figure}
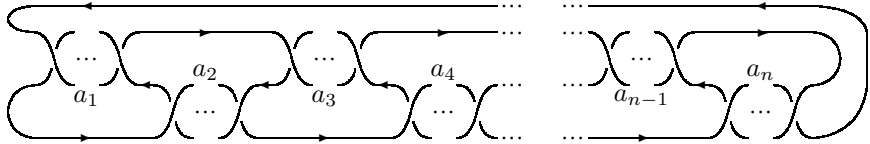

Here, we are only interested in knots, so let us regard the $2$-bridge knots
as those $S(\alpha,\beta)$ with presentations $C(a_1,\dots,a_n)$ where both
$n$ and each $a_i$ are even (and nonzero!). We shall refer to such a
presentation as an \emph{even} (Conway) presentation. In fact, all $2$-bridge
knots $S(\alpha,\beta)$ are necessarily \emph{prime}, and it turns out that
they have odd $\alpha$. We assign an \emph{orientation} to an even
presentation $C(a)$ by giving that of the upper arc, as depicted in
Figure~\ref{figure:EvenConwayPresentationForOrientedTwoBridgeKnots}.

Next, note the equality: $C(a_1,\dots,a_n)=C(a_1,\dots, a_n \pm 1,\mp 1)$.
More generally, if none of the first and last components of sequences $a$ and
$b$ are of unit magnitude, the oriented $2$-bridge links with presentations
$C(a)$ and $C(b)$ are equal if and only if $a$ and $b$ are of the same length
$n$ and $b=a$ or $b=(-)^{n+1}a^r$, where by $a^r$ we intend the sequence
$a_n,\dots,a_1$. This fact means that an even presentation $C(a)$,
corresponding to some $2$-bridge knot type $S(\alpha,\beta)$ is in fact a
\emph{unique} even presentation for $S(\alpha,\beta)$, modulo only the
consideration that $C(a)^*=C(-a)$. If we consider $C(a)$ as oriented, then by
inspection, the inverse is: $-C(a)=C(-a^r)$. However $2$-bridge links are
known to be invertible, so $C(a)=C(-a^r)$ and $C(-a)=C(a^r)$; and if $a$ is
palindromic, that is if $a=a^r$, then $C(a)$ is achiral and in fact
$C(a)=C(a^r)=C(-a)=C(-a^r)$.


\section{DT codes for Conway presentations}

We now describe how to construct a Dowker--Thistlethwaite code corresponding
to an even Conway presentation $C(a_1,\dots,a_n)$. A fluent description of the
algorithm used to construct such codes is available in
\cite{HosteThistlethwaiteWeeks:1998}. The algorithm begins by walking once
around an oriented diagram, numbering crossings from $1$ onwards, and
multiplying by $-1$ any even numbers written down whilst travelling
\emph{over} a crossing. So, each crossing of an $m$-crossing presentation
bears one unsigned odd and one signed even number, drawn from the set
$\{1,\dots,2m\}$. The DT code is then formed as the sequence of the even
numbers corresponding to ordered odd crossing numbers $1,3,\dots,2m-1$.

To apply this procedure to an even Conway presentation $C(a_1,\dots,a_n)$, we
initially number the crossings increasing from the upper left of
Figure~\ref{figure:EvenConwayPresentationForOrientedTwoBridgeKnots}. Then, for
each $i=1,\dots,n$, let $m_i$ be the (left-to-right) sequence associated with
the leftwards-travelling middle strand and $o_i$ the sequence associated with
the rightwards-travelling outer strand (which alternates with $i$ between the
upper and lower strands). Then, the DT code associated with $C(a)$ describes
the pairings $(m_i,o_i)$ for $i=1,\dots,n$. Setting $A_i:=|a_i|$ and
$s_i:=\mathrm{sign}(a_i)$ for $i=1,\dots,n$, by inspection, the components of
the sequences are:
\begin{eqnarray*}
  \begin{array}{rcl}
    (m_i)_k
    & = &
    (M_i + 1 - k) s_i^k
    \\
    (o_i)_k
    & = &
    (O_i + k) s_i^{k+1}
  \end{array}
  \qquad
  \mathrm{for}
  \;\;
  k=1,\dots,A_i,
\end{eqnarray*}
where the offsets $M_i$ and $O_i$ are:
\begin{eqnarray*}
  M_i
  & = &
  O_{n-1}
  +
  {\textstyle \sum_{j=i}^n} A_j
  \\
  O_i
  & = &
  \left\{
    \begin{array}{ll}
      \sum_{j=1}^{i-2,\mathrm{odd}} A_j,
      &
      \mathrm{odd}
      \;
      i
      \\
      M_1
      +
      \sum_{j=2}^{i-2,\mathrm{even}} A_j,
      &
      \mathrm{even}
      \;
      i.
    \end{array}
  \right.
\end{eqnarray*}

For example, for $C(+2,-2)$, which corresponds to the (negative, left-handed)
trefoil knot $3^A_1$, we obtain the DT code $(6,-8,2,-4)$; note that the DT
code may be regarded as corresponding to either chirality of $3^A_1$. In
contrast the palindromic presentation $C(+2,+2)$ is the achiral $4^A_1$ with
DT code $(6,8,2,4)$.


\section{Enumerating even Conway presentations}

For even positive integers $m$ and $n$, let $U^{m,n}$ be the set of even
Conway presentations $C(a_1,\dots,a_n)$ of $m=\sum_{i=1}^n |a_i|$ crossings
(so that $n\leqslant m/2$), restricted such that $U^{m,n}$ contains only one
element of each knot type equivalence class $\{C(a),C(a^r),C(-a),C(-a^r)\}$.
Denote by $TK^{m,n}$ the size of $U^{m,n}$ and set $TK_c^{m,n}$ as the number
of $c$-crossing knots (not presentations) in $U^{m,n}$. Then we have $TK^{m,n}
= \sum_c TK_c^{m,n}$, where the index $c$ runs over some finite set of
crossing numbers depending on $m$ and $n$, and $TK_c = \sum_{m,n} TK_c^{m,n}$,
where the indices $m,n$ run over some finite set depending on $c$.

An efficient procedure (which we have implemented in \textsc{Mathematica}) to
compile the sequences underlying $U^{m,n}$ is to:
\begin{enumerate}
  \item
    find all ordered sets $\{a_1,\dots,a_n\}$ of even \emph{positive} integers
    $a_i$ such that $\sum_{i=1}^n a_i=m$, then
  \item
    find all permutations of each such set, and
  \item
    for each such permutation, find all $2^n$ variants obtained by
    multiplying each component by $\pm 1$, and finally
  \item
    identify all cliques of the form $\{a,a^r,-a,-a^r\}$, and eliminate
    all but one element from each (for good measure, select one with
    $a_1>0$).
\end{enumerate}

Here we are interested in even presentations $C(a)$ which are prime knots of
at most $16$ crossings, and by inspection (see below) we find that for all
even $m$ between $4$ and $28$ there exist some even $n\leqslant 14$ such that
$U^{m,n}$ contains some such $C(a)$. Table~\ref{table:TKmn} lists $TK^{m,n}$
for various $m$ and $n$ of interest, together with their sums $TK^m \triangleq
\sum_n TK^{m,n}$, where the sum runs over all even positive $n$ such that $2n
\leqslant m$.

\renewcommand{\arraystretch}{1.5}
\renewcommand{\tabcolsep}{4pt}

\begin{table}[ht]
  \tiny
  \begin{centering}
  \begin{tabular}{r|*{7}{r}|r}
    \multicolumn{1}{c}{} & \multicolumn{7}{c}{$n$} & \\
    $m$ & 2&   4&     6&  8& 10&12&     14& $TK^m$ \\
    \hline
      4 & 2&   .&     .&  .&  .& .& .&      2 \\
      6 & 2&   .&     .&  .&  .& .& .&      2 \\
      8 & 4&   6&     .&  .&  .& .& .&     10 \\
     10 & 4&  16&     .&  .&  .& .& .&     20 \\
     12 & 6&  44&    20&  .&  .& .& .&     70 \\
     14 & 6&  80&    96&  .&  .& .& .&    182 \\
     16 & 8& 146&   348& 72&  .& .& .&    574 \\
    \cline{2-2}
     18&\multicolumn{1}{r|}{  8}&   224&   896&    512&  .& .& .&   1640 \\
    \cline{3-3}
     20& 10&\multicolumn{1}{r|}{   344}&  2040&   2336&    272& .& .&   5002 \\
    \cline{4-4}
     22& 10&   480&\multicolumn{1}{r|}{   4032}&   7680&   2560& .& .&  14762 \\
    \cline{5-5}
     24& 12&  670&  7432&\multicolumn{1}{r|}{  21200}&  14160&  1056& .&  44530 \\
    \cline{6-6}
     26& 12&  880&  12672& 50688&\multicolumn{1}{r|}{  56320}&  12288&  .&  132860 \\
    \cline{7-7}
     28& 14&  1156& 20652& 109984& 183280&\multicolumn{1}{r|}{ 80064}&  4160& 399310 \\
    \hline
     total& 98&  4046& 48188& 192472& 256592& 93408&  4160& 598964
  \end{tabular}
  \caption{%
    The numbers $TK^{m,n}$ of $2$-bridge knot types defined by $m$-crossing even presentations
    $C(a_1,\dots,a_n)$, together with their sums $TK^m$.
    Above the zigzag are $U^{m,n}$ which contribute to our census
    of $2$-bridge knots of up to $16$ crossings.
  }%
  \label{table:TKmn}
  \end{centering}
\end{table}

We use \textsc{Mathematica} to compute DT codes for all elements of $U^{m,n}$
for the various $m$ and $n$ described in Table~\ref{table:TKmn}. These DT
codes are then fed to the \texttt{Locate} function of \textsc{Knotscape} for
identification as alternating prime knots within the HTW tables. For the vast
majority of the presentations with $m>16$, \textsc{Knotscape} reports that the
knot is of more than $16$ crossings, and instead of returning an
identification, it returns a reduced DT code and associated crossing number.

Inspection of the results of a range of searches of various $U^{m,n}$
indicates that $T_c^{m,n}$ is nonzero exactly for $m-n+1\leqslant c \leqslant
m$. This shows that for a given $m$, presentations with higher values of $n$
have greater opportunities for reduction. For our purposes, we need not prove
that $T_c^{m,n}=0$ for $c<m-n+1$ (see below), however it appears that a closer
reading of~\cite{ErnstSumners:1987} should indicate its origin. Demanding that
the crossing numbers $c$ be at most $16$ thus means demanding $m-n+1 \leqslant
16$. As both $m$ and $n$ are even and $2n\leqslant m$, we then have
$m\leqslant 28$ and $n\geqslant m-14$; that is, only choices of $m$ and $n$
above the zigzag line in Table~\ref{table:TKmn} contribute to our collection
of $2$-bridge knots of up to $16$ crossings.


\section{The results}

In this manner, we identify as $2$-bridge knots a total of $5546$ alternating
prime knots within the HTW tables. That is, corresponding to each such knot
type, as defined by its DT code, we have an oriented $2$-bridge knot defined
by an even Conway presentation. As mentioned above, the \texttt{Locate}
function identifies a given DT code as corresponding to some knot \emph{type}
in the tables, and it can do no more as a DT code corresponds to both a knot
and its reflection (although it does determine their orientations). The $5546$
knots we identify become $11007$ when we count knots and their reflections as
distinct; that is, there are $85$ achiral knots which correspond precisely to
the even presentations $C(a)$ with palindromic sequences $a$.

Tables~\ref{table:TKcmn} and~\ref{table:TKcmandTKc} describe $TK^{m,n}_c$,
$TK_c^m \triangleq \sum_n TK_c^{m,n}$ and $TK_c = \sum_m TK_c^m$ for $3
\leqslant c \leqslant 16$ for the $U^{m,n}$ mentioned in
Table~\ref{table:TKmn}. The numbers in the last row of Table~\ref{table:TKcmn}
sum to $598,964$, the number of $2$-bridge knot types listed in
Table~\ref{table:TKmn}. For $m\leqslant 16$ (only), the final column of
Table~\ref{table:TKcmandTKc} reconstructs $TK^m=\sum_c TK_c^m$, in agreement
with the final column of Table~\ref{table:TKmn}. Similarly, for $m\leqslant
16$ (only), the final row of Table~\ref{table:TKcmn} reconstructs
$TK_c=\sum_{m,n} TK_c^{m,n}$, in agreement with the final row of
Table~\ref{table:TKc}. (For $m>16$, we do not have enough data to recover
these numbers.)

\renewcommand{\arraystretch}{1.5}
\renewcommand{\tabcolsep}{1.5pt}

\begin{table}[htbp]
  \tiny
  \hspace{-40pt}
  \begin{tabular}{r@{\hspace{4pt}}|@{\hspace{4pt}}*{26}{r}}
                 & \multicolumn{26}{c}{$c$} \\
     $m, n$& 3& 4& 5& 6& 7&  8&  9& 10& 11&  12&  13&  14&   15&   16&   17&    18&    19&    20&    21&    22&    23&     24&     25&    26&    27&   28 \\
   \hline
     4, 2& 1& 1&     . &     . &     . &     .  &     .  &     .  &     .  &     .   &     .   &     .   &     .    &     .    &     .    &     .     &     .     &     .     &     .     &     .     &     .      &     .      &     .     &     .     &     .     &     .     \\
     6, 2&     . &     . & 1& 1&     . &     .  &     .  &     .  &     .  &     .   &     .   &     .   &     .    &     .    &     .    &     .     &     .     &     .     &     .     &     .     &     .      &     .      &     .     &     .     &     .     &     .     \\
     8, 2&     . &     . &     . &     . & 2&  2&     .  &     .  &     .  &     .   &     .   &     .   &     .    &     .    &     .    &     .     &     .     &     .     &     .     &     .     &     .      &     .      &     .     &     .     &     .     &     .     \\
     8, 4&     . &     . & 1& 2& 2&  1&     .  &     .  &     .  &     .   &     .   &     .   &     .    &     .    &     .    &     .     &     .     &     .     &     .     &     .     &     .      &     .      &     .     &     .     &     .     &     .     \\
    10, 2&     . &     . &     . &     . &     . &     .  &  2&  2&     .  &     .   &     .   &     .   &     .    &     .    &     .    &     .     &     .     &     .     &     .     &     .     &     .      &     .      &     .     &     .     &     .     &     .     \\
    10, 4&     . &     . &     . &     . & 2&  6&  6&  2&     .  &     .   &     .   &     .   &     .    &     .    &     .    &     .     &     .     &     .     &     .     &     .     &     .      &     .      &     .     &     .     &     .     &     .     \\
    12, 2&     . &     . &     . &     . &     . &     .  &     .  &     .  &  3&   3&     .   &     .   &     .    &     .    &     .    &     .     &     .     &     .     &     .     &     .     &     .      &     .      &     .     &     .     &     .     &     .     \\
    12, 4&     . &     . &     . &     . &     . &     .  &  6& 16& 16&   6&     .   &     .   &     .    &     .    &     .    &     .     &     .     &     .     &     .     &     .     &     .      &     .      &     .     &     .     &     .     &     .     \\
    12, 6&     . &     . &     . &     . & 1&  3&  6&  6&  3&   1&     .   &     .   &     .    &     .    &     .    &     .     &     .     &     .     &     .     &     .     &     .      &     .      &     .     &     .     &     .     &     .     \\
    14, 2&     . &     . &     . &     . &     . &     .  &     .  &     .  &     .  &     .   &   3&   3&     .    &     .    &     .    &     .     &     .     &     .     &     .     &     .     &     .      &     .      &     .     &     .     &     .     &     .     \\
    14, 4&     . &     . &     . &     . &     . &     .  &     .  &     .  & 10&  30&  30&  10&     .    &     .    &     .    &     .     &     .     &     .     &     .     &     .     &     .      &     .      &     .     &     .     &     .     &     .     \\
    14, 6&     . &     . &     . &     . &     . &     .  &  3& 15& 30&  30&  15&   3&     .    &     .    &     .    &     .     &     .     &     .     &     .     &     .     &     .      &     .      &     .     &     .     &     .     &     .     \\
    16, 2&     . &     . &     . &     . &     . &     .  &     .  &     .  &     .  &     .   &     .   &     .   &    4&    4&     .    &     .     &     .     &     .     &     .     &     .     &     .      &     .      &     .     &     .     &     .     &     .     \\
    16, 4&     . &     . &     . &     . &     . &     .  &     .  &     .  &     .  &     .   &  19&  54&   54&   19&     .    &     .     &     .     &     .     &     .     &     .     &     .      &     .      &     .     &     .     &     .     &     .     \\
    16, 6&     . &     . &     . &     . &     . &     .  &     .  &     .  & 12&  54& 108& 108&   54&   12&     .    &     .     &     .     &     .     &     .     &     .     &     .      &     .      &     .     &     .     &     .     &     .     \\
    16, 8&     . &     . &     . &     . &     . &     .  &  1&  4& 12&  19&  19&  12&    4&    1&     .    &     .     &     .     &     .     &     .     &     .     &     .      &     .      &     .     &     .     &     .     &     .     \\
    18, 2&     . &     . &     . &     . &     . &     .  &     .  &     .  &     .  &     .   &     .   &     .   &     .    &     .    &    4&     4&     .     &     .     &     .     &     .     &     .      &     .      &     .     &     .     &     .     &     .     \\
    18, 4&     . &     . &     . &     . &     . &     .  &     .  &     .  &     .  &     .   &     .   &     .   &   28&   84&   84&    28&     .     &     .     &     .     &     .     &     .      &     .      &     .     &     .     &     .     &     .     \\
    18, 6&     . &     . &     . &     . &     . &     .  &     .  &     .  &     .  &     .   &  28& 140&  280&  280&  140&    28&     .     &     .     &     .     &     .     &     .      &     .      &     .     &     .     &     .     &     .     \\
    18, 8&     . &     . &     . &     . &     . &     .  &     .  &     .  &  4&  28&  84& 140&  140&   84&   28&     4&     .     &     .     &     .     &     .     &     .      &     .      &     .     &     .     &     .     &     .     \\
    20, 2&     . &     . &     . &     . &     . &     .  &     .  &     .  &     .  &     .   &     .   &     .   &     .    &     .    &     .    &     .     &     5&     5&     .     &     .     &     .      &     .      &     .     &     .     &     .     &     .     \\
    20, 4&     . &     . &     . &     . &     . &     .  &     .  &     .  &     .  &     .   &     .   &     .   &     .    &     .    &   44&   128&   128&    44&     .     &     .     &     .      &     .      &     .     &     .     &     .     &     .     \\
    20, 6&     . &     . &     . &     . &     . &     .  &     .  &     .  &     .  &     .   &     .   &     .   &   66&  318&  636&   636&   318&    66&     .     &     .     &     .      &     .      &     .     &     .     &     .     &     .     \\
    20, 8&     . &     . &     . &     . &     . &     .  &     .  &     .  &     .  &     .   &  20& 128&  384&  636&  636&   384&   128&    20&     .     &     .     &     .      &     .      &     .     &     .     &     .     &     .     \\
    20,10&     . &     . &     . &     . &     . &     .  &     .  &     .  &  1&   5&  20&  44&   66&   66&   44&    20&     5&     1&     .     &     .     &     .      &     .      &     .     &     .     &     .     &     .     \\
    22, 2&     . &     . &     . &     . &     . &     .  &     .  &     .  &     .  &     .   &     .   &     .   &     .    &     .    &     .    &     .     &     .     &     .     &     5&     5&     .      &     .      &     .     &     .     &     .     &     .     \\
    22, 4&     . &     . &     . &     . &     . &     .  &     .  &     .  &     .  &     .   &     .   &     .   &     .    &     .    &     .    &     .     &    60&   180&   180&    60&     .      &     .      &     .     &     .     &     .     &     .     \\
    22, 6&     . &     . &     . &     . &     . &     .  &     .  &     .  &     .  &     .   &     .   &     .   &     .    &     .    &  126&   630&  1260&  1260&   630&   126&     .      &     .      &     .     &     .     &     .     &     .     \\
    22, 8&     . &     . &     . &     . &     . &     .  &     .  &     .  &     .  &     .   &     .   &     .   &   60&  420& 1260&  2100&  2100&  1260&   420&    60&     .      &     .      &     .     &     .     &     .     &     .     \\
    22,10&     . &     . &     . &     . &     . &     .  &     .  &     .  &     .  &     .   &   5&  45&  180&  420&  630&   630&   420&   180&    45&     5&     .      &     .      &     .     &     .     &     .     &     .     \\
    24, 2&     . &     . &     . &     . &     . &     .  &     .  &     .  &     .  &     .   &     .   &     .   &     .    &     .    &     .    &     .     &     .     &     .     &     .     &     .     &      6&      6&     .     &     .     &     .     &     .     \\
    24, 4&     . &     . &     . &     . &     . &     .  &     .  &     .  &     .  &     .   &     .   &     .   &     .    &     .    &     .    &     .     &     .     &     .     &    85&   250&    250&     85&     .     &     .     &     .     &     .     \\
    24, 6&     . &     . &     . &     . &     . &     .  &     .  &     .  &     .  &     .   &     .   &     .   &     .    &     .    &     .    &     .     &   236&  1160&  2320&  2320&   1160&    236&     .     &     .     &     .     &     .     \\
    24, 8&     . &     . &     . &     . &     . &     .  &     .  &     .  &     .  &     .   &     .   &     .   &     .    &     .    &  170&  1160&  3480&  5790&  5790&  3480&   1160&    170&     .     &     .     &     .     &     .     \\
    24,10&     . &     . &     . &     . &     . &     .  &     .  &     .  &     .  &     .   &     .   &     .   &   30&  250& 1000&  2320&  3480&  3480&  2320&  1000&    250&     30&     .     &     .     &     .     &     .     \\
    24,12&     . &     . &     . &     . &     . &     .  &     .  &     .  &     .  &     .   &   1&   6&   30&   85&  170&   236&   236&   170&    85&    30&      6&      1&     .     &     .     &     .     &     .     \\
    26, 2&     . &     . &     . &     . &     . &     .  &     .  &     .  &     .  &     .   &     .   &     .   &     .    &     .    &     .    &     .     &     .     &     .     &     .     &     .     &     .      &     .      &     6&     6&     .     &     .     \\
    26, 4&     . &     . &     . &     . &     . &     .  &     .  &     .  &     .  &     .   &     .   &     .   &     .    &     .    &     .    &     .     &     .     &     .     &     .     &     .     &    110&    330&   330&   110&     .     &     .     \\
    26, 6&     . &     . &     . &     . &     . &     .  &     .  &     .  &     .  &     .   &     .   &     .   &     .    &     .    &     .    &     .     &     .     &     .     &   396&  1980&   3960&   3960&  1980&   396&     .     &     .     \\
    26, 8&     . &     . &     . &     . &     . &     .  &     .  &     .  &     .  &     .   &     .   &     .   &     .    &     .    &     .    &     .     &   396&  2772&  8316& 13860&  13860&   8316&  2772&   396&     .     &     .     \\
    26,10&     . &     . &     . &     . &     . &     .  &     .  &     .  &     .  &     .   &     .   &     .   &     .    &     .    &  110&   990&  3960&  9240& 13860& 13860&   9240&   3960&   990&   110&     .     &     .     \\
    26,12&     . &     . &     . &     . &     . &     .  &     .  &     .  &     .  &     .   &     .   &     .   &    6&   66&  330&   990&  1980&  2772&  2772&  1980&    990&    330&    66&     6&     .     &     .     \\
    28, 2&     . &     . &     . &     . &     . &     .  &     .  &     .  &     .  &     .   &     .   &     .   &     .    &     .    &     .    &     .     &     .     &     .     &     .     &     .     &     .      &     .      &     .     &     .     &     7&    7 \\
    28, 4&     . &     . &     . &     . &     . &     .  &     .  &     .  &     .  &     .   &     .   &     .   &     .    &     .    &     .    &     .     &     .     &     .     &     .     &     .     &     .      &     .      &   146&   432&   432&  146 \\
    28, 6&     . &     . &     . &     . &     . &     .  &     .  &     .  &     .  &     .   &     .   &     .   &     .    &     .    &     .    &     .     &     .     &     .     &     .     &     .     &    651&   3225&  6450&  6450&  3225&  651 \\
    28, 8&     . &     . &     . &     . &     . &     .  &     .  &     .  &     .  &     .   &     .   &     .   &     .    &     .    &     .    &     .     &     .     &     .     &   868&  6016&  18048&  30060& 30060& 18048&  6016&  868 \\
    28,10&     . &     . &     . &     . &     . &     .  &     .  &     .  &     .  &     .   &     .   &     .   &     .    &     .    &     .    &     .     &   365&  3225& 12900& 30060&  45090&  45090& 30060& 12900&  3225&  365 \\
    28,12&     . &     . &     . &     . &     . &     .  &     .  &     .  &     .  &     .   &     .   &     .   &     .    &     .    &   42&   432&  2160&  6450& 12900& 18048&  18048&  12900&  6450&  2160&   432&   42 \\
    28,14&     . &     . &     . &     . &     . &     .  &     .  &     .  &     .  &     .   &     .   &     .   &    1&    7&   42&   146&   365&   651&   868&   868&    651&    365&   146&    42&     7&    1 \\
   \hline
    totals& 1& 1& 2& 3& 7& 12& 24& 45& 91& 176& 352& 693& 1387& 2752& 5496& 10866& 21082& 38726& 64760& 94008& 113480& 109064& 79456& 41056& 13344& 2080
  \end{tabular}
  \caption{%
    The numbers $TK_c^{m,n}$ of $c$-crossing alternating prime knots identified as
    $m$-crossing $2$-bridge knot presentations $C(a_1,\dots,a_n)$ within
    the $U^{m,n}$ of Table~\ref{table:TKmn}.
  }%
  \label{table:TKcmn}
\end{table}

\renewcommand{\arraystretch}{1.5}
\renewcommand{\tabcolsep}{4pt}

\begin{table}[htbp]
  \tiny
  \begin{centering}
  \begin{tabular}{r|*{14}{r}|r}
    \multicolumn{1}{c}{} & \multicolumn{14}{c}{$c$} \\
      $m$&     3&     4&     5&     6&     7&     8&     9&    10&    11&    12&    13&    14&    15&    16& totals \\
    \hline
        4&     1&     1&     .&     .&     .&     .&     .&     .&     .&     .&     .&     .&     .&     .&              2 \\
        6&     .&     .&     1&     1&     .&     .&     .&     .&     .&     .&     .&     .&     .&     .&              2 \\
        8&     .&     .&     1&     2&     4&     3&     .&     .&     .&     .&     .&     .&     .&     .&             10 \\
       10&     .&     .&     .&     .&     2&     6&     8&     4&     .&     .&     .&     .&     .&     .&             20 \\
       12&     .&     .&     .&     .&     1&     3&    12&    22&    22&    10&     .&     .&     .&     .&             70 \\
       14&     .&     .&     .&     .&     .&     .&     3&    15&    40&    60&    48&    16&     .&     .&            182 \\
       16&     .&     .&     .&     .&     .&     .&     1&     4&    24&    73&   146&   174&   116&    36&            574 \\
       18&     .&     .&     .&     .&     .&     .&     .&     .&     4&    28&   112&   280&   448&   448&           1320 \\
       20&     .&     .&     .&     .&     .&     .&     .&     .&     1&     5&    40&   172&   516&  1020&           1754 \\
       22&     .&     .&     .&     .&     .&     .&     .&     .&     .&     .&     5&    45&   240&   840&           1130 \\
       24&     .&     .&     .&     .&     .&     .&     .&     .&     .&     .&     1&     6&    60&   335&            402 \\
       26&     .&     .&     .&     .&     .&     .&     .&     .&     .&     .&     .&     .&     6&    66&             72 \\
       28&     .&     .&     .&     .&     .&     .&     .&     .&     .&     .&     .&     .&     1&     7&              8 \\
    \hline
   $TK_c$&     1&     1&     2&     3&     7&    12&    24&    45&    91&   176&   352&   693&  1387&   2752&          5546
  \end{tabular}
  \caption{%
    Observed $TK_c^m = \sum_n TK_c^{m,n}$, and $TK_c = \sum_m TK_c^m$.
  }%
  \label{table:TKcmandTKc}
  \end{centering}
\end{table}

Importantly, our list of $5546$ knots is \emph{complete}, as our observed
$TK_c$ are in agreement with those of Table~\ref{table:TKc}. (For this reason,
we need not prove the experimental observation that $T_c^{m,n}=0$ for
$c<m-n+1$.) This fact, together with the symmetries of
Table~\ref{table:TKcmn}, indicates that the \texttt{Locate} function is
successful in reducing to minimal presentations all $598,964$ of the
$2$-bridge knot types with even Conway presentations of up to $28$ crossings.
This success is a remarkable credit to the ability of the \texttt{Locate}
function to reduce an input $m$-term DT code corresponding to a $c$-crossing
alternating knot to a $c$-term DT code (where we intend $c\leqslant 16$). The
\texttt{Locate} function (of course!) cannot be guaranteed to always succeed
in this if $m>16$.

To illustrate the results, Table~\ref{table:2bridgeknotdata} lists the
$2$-bridge knots of up to $8$ crossings, together with the sequences $a$ of
their even Conway presentations $C(a)$, and associated DT codes.
Tables~\ref{table:Thetwobridgeknotsof3to14crossings}--\ref{table:Thetwobridgeknotsof16crossings2}
then list the indices of all the $2$-bridge knots of between $3$ and $16$
crossings. Sequences $C(a)$ for these knots are available on request from the
author.

\renewcommand{\arraystretch}{1.5}
\renewcommand{\tabcolsep}{4pt}

\begin{table}[htbp]
  \tiny
  \begin{centering}
  \begin{tabular}{r|r@{\hspace{4pt}}|@{\hspace{4pt}}*{6}{r}@{\hspace{4pt}}|@{\hspace{4pt}}*{12}{r}}
    \multicolumn{1}{c}{knot} & \multicolumn{1}{c}{$m,n$} & \multicolumn{6}{c}{$a$} & \multicolumn{12}{c}{DT~code} \\
    \hline
     $3^A_{ 1}$& 4,2& 2& -2& & & & & 6& -8& 2& -4& & & & & & & & \\
     $4^A_{ 1}$& 4,2& 2& 2& & & & & 6& 8& 2& 4& & & & & & & & \\
     $5^A_{ 1}$& 6,2& 4& -2& & & & & 10& 8& -12& 4& 2& -6& & & & & & \\
     $5^A_{ 2}$& 8,4& 2& -2& 2& -2& & & 12& 8& -16& 4& -14& 2& -10& -6& & & & \\
     $6^A_{ 1}$& 8,4& 2& -2& -2& 2& & & 12& -8& 16& -4& -14& 2& -10& 6& & & & \\
     $6^A_{ 2}$& 8,4& 2& 2& -2& 2& & & 12& -8& 16& -4& 14& 2& 10& 6& & & & \\
     $6^A_{ 3}$& 6,2& 4& 2& & & & & 10& 8& 12& 4& 2& 6& & & & & & \\
     $7^A_{ 1}$& 8,4& 2& 2& -2& -2& & & 12& -8& -16& -4& 14& 2& 10& -6& & & & \\
     $7^A_{ 2}$& 8,4& 2& 2& 2& -2& & & 12& 8& -16& 4& 14& 2& 10& -6& & & & \\
     $7^A_{ 3}$& 10,4& 2& -2& 4& -2& & & 16& 12& 10& -20& 6& 4& -18& 2& -14& -8& & \\
     $7^A_{ 4}$& 8,2& 6& -2& & & & & 14& 12& 10& -16& 6& 4& 2& -8& & & & \\
     $7^A_{ 5}$& 10,4& 4& -2& 2& -2& & & 16& 14& 10& -20& 6& -18& 4& 2& -12& -8& & \\
     $7^A_{ 6}$& 8,2& 4& -4& & & & & 12& 10& -16& -14& 4& 2& -8& -6& & & & \\
     $7^A_{ 7}$& 12,6& 2& -2& 2& -2& 2& -2 & 18& 14& 10& -24& 6& -22& 4& -20& 2& -16& -12& -8 \\
     $8^A_{ 1}$& 10,4& 2& 2& -4& 2& & & 16& -12& -10& 20& -6& -4& 18& 2& 14& 8& & \\
     $8^A_{ 4}$& 10,4& 2& -2& -4& 2& & & 16& -12& -10& 20& -6& -4& -18& 2& -14& 8& & \\
     $8^A_{ 5}$& 8,4& 2& 2& 2& 2& & & 12& 8& 16& 4& 14& 2& 10& 6& & & & \\
     $8^A_{ 6}$& 12,6& 2& -2& 2& -2& -2& 2 & 18& 14& -10& 24& -6& -22& 4& -20& 2& -16& -12& 8 \\
     $8^A_{ 7}$& 10,4& 4& -2& -2& 2& & & 16& 14& -10& 20& -6& -18& 4& 2& -12& 8& & \\
     $8^A_{ 8}$& 12,6& 2& 2& -2& 2& -2& 2 & 18& -14& -10& 24& -6& 22& -4& 20& 2& 16& 12& 8 \\
     $8^A_{ 9}$& 10,4& 4& -2& 2& 2& & & 16& 14& 10& 20& 6& -18& 4& 2& -12& 8& & \\
     $8^A_{10}$& 10,4& 2& 4& -2& 2& & & 14& -8& 20& -4& 18& 16& 2& 12& 10& 6& & \\
     $8^A_{11}$& 8,2& 6& 2& & & & & 14& 12& 10& 16& 6& 4& 2& 8& & & & \\
     $8^A_{16}$& 12,6& 2& -2& 2& 2& -2& 2 & 18& 14& -10& 24& -6& 22& 4& -20& 2& -16& 12& 8 \\
     $8^A_{17}$& 10,4& 4& 2& -2& 2& & & 16& 14& -10& 20& -6& 18& 4& 2& 12& 8& & \\
     $8^A_{18}$& 8,2& 4& 4& & & & & 12& 10& 16& 14& 4& 2& 8& 6& & & & \\
  \end{tabular}
  \caption{%
    The $2$-bridge knots of up to $8$ crossings, with even Conway presentations $C(a) \in
    U^{m,n}$, and associated DT codes.
  }%
  \label{table:2bridgeknotdata}
  \end{centering}
\end{table}

\begin{table}[htbp]
  \tiny
  \hspace{-50pt}
  \begin{tabular}{p{150mm}}
    \multicolumn{1}{c}{\normalsize $ 8$ crossings} \\[1.5mm]
    1, 4, 5, 6, 7, 8, 9, 10, 11, 16, 17, 18
    \\[1.5mm]
    \multicolumn{1}{c}{\normalsize $ 9$ crossings} \\[1.5mm]
    3, 8, 10, 12, 13, 14, 15, 16, 17, 19, 20, 21, 22, 23, 24, 26, 27, 33, 34, 35, 36, 38, 39, 41
    \\[1.5mm]
    \multicolumn{1}{c}{\normalsize $10$ crossings} \\[1.5mm]
    5, 19, 23, 25, 26, 29, 30, 31, 32, 33, 34, 35, 43, 44, 49, 52, 53, 54, 55, 56, 57, 58, 59, 60, 61, 63, 64, 65, 68, 69, 70, 71, 74, 75, 107, 108, 109, 110, 111, 112, 113, 114, 115, 116, 117
    \\[1.5mm]
    \multicolumn{1}{c}{\normalsize $11$ crossings} \\[1.5mm]
    13, 59, 65, 75, 77, 84, 85, 89, 90, 91, 93, 95, 96, 98, 110, 111, 117, 119, 120, 121, 140, 144, 145, 154, 159, 166, 174, 175, 176, 177, 178, 179, 180, 182, 183, 184, 185, 186, 188, 190, 191, 192, 193, 195, 203, 204, 205, 206, 207, 208, 210, 211, 220, 224, 225, 226, 229, 230, 234, 235, 236, 238, 242, 243, 246, 247, 306, 307, 308, 309, 310, 311, 333, 334, 335, 336, 337, 339, 341, 342, 343, 355, 356, 357, 358, 359, 360, 363, 364, 365, 367
    \\[1.5mm]
    \multicolumn{1}{c}{\normalsize $12$ crossings} \\[1.5mm]
    38, 169, 197, 204, 206, 221, 226, 239, 241, 243, 247, 251, 254, 255, 257, 259, 300, 302, 303, 306, 307, 330, 378, 379, 380, 384, 385, 406, 425, 437, 447, 454, 471, 477, 482, 497, 498, 499, 500, 501, 502, 506, 508, 510, 511, 512, 514, 517, 518, 519, 520, 521, 522, 528, 532, 533, 534, 535, 536, 537, 538, 539, 540, 541, 545, 549, 550, 551, 552, 579, 580, 581, 582, 583, 584, 585, 595, 596, 597, 600, 601, 643, 644, 649, 650, 651, 652, 682, 684, 690, 691, 713, 714, 715, 716, 717, 718, 720, 721, 722, 723, 724, 726, 727, 728, 729, 731, 732, 733, 736, 738, 740, 743, 744, 745, 758, 759, 760, 761, 762, 763, 764, 773, 774, 775, 791, 792, 796, 797, 802, 803, 1023, 1024, 1029, 1030, 1033, 1034, 1039, 1040, 1125, 1126, 1127, 1128, 1129, 1130, 1131, 1132, 1133, 1134, 1135, 1136, 1138, 1139, 1140, 1145, 1146, 1148, 1149, 1157, 1158, 1159, 1161, 1162, 1163, 1165, 1166, 1273, 1274, 1275, 1276, 1277, 1278, 1279, 1281, 1282, 1287
    \\[1.5mm]
    \multicolumn{1}{c}{\normalsize $13$ crossings} \\[1.5mm]
    110, 640, 667, 711, 718, 776, 778, 829, 830, 847, 848, 853, 857, 864, 869, 871, 876, 880, 1003, 1004, 1038, 1046, 1048, 1050, 1241, 1242, 1284, 1293, 1295, 1296, 1461, 1474, 1475, 1476, 1540, 1573, 1585, 1608, 1619, 1640, 1647, 1739, 1741, 1743, 1749, 1762, 1765, 1766, 1768, 1772, 1773, 1774, 1777, 1779, 1780, 1783, 1784, 1785, 1789, 1791, 1793, 1796, 1797, 1799, 1808, 1812, 1858, 1859, 1860, 1861, 1866, 1867, 1868, 1869, 1876, 1877, 1879, 1880, 1881, 1882, 1883, 1884, 1888, 1889, 1892, 1904, 1920, 1922, 1925, 1928, 1930, 1931, 1934, 1937, 1938, 1942, 1943, 1944, 1946, 1948, 1951, 1975, 1977, 2047, 2048, 2049, 2050, 2051, 2067, 2070, 2072, 2073, 2099, 2100, 2112, 2114, 2115, 2116, 2120, 2121, 2123, 2124, 2125, 2126, 2137, 2145, 2146, 2280, 2281, 2288, 2289, 2290, 2291, 2318, 2319, 2445, 2458, 2464, 2465, 2475, 2478, 2479, 2574, 2577, 2584, 2675, 2676, 2677, 2678, 2679, 2680, 2683, 2684, 2685, 2686, 2687, 2688, 2689, 2690, 2692, 2693, 2696, 2697, 2698, 2699, 2700, 2701, 2702, 2703,
    2704, 2705, 2709, 2716, 2717, 2718, 2719, 2721, 2722, 2723, 2724, 2725, 2726, 2734, 2751, 2752, 2753, 2754, 2755, 2756, 2757, 2758, 2759, 2760, 2761, 2762, 2763, 2776, 2777, 2778, 2786, 2791, 2818, 2819, 2820, 2823, 2824, 2825, 2826, 2827, 2828, 2829, 2830, 2831, 2833, 2834, 2835, 2836, 2842, 2843, 2854, 2855, 2858, 2859, 2863, 2865, 2874, 2875, 2934, 2935, 2945, 2946, 2947, 2948, 2949, 2950, 2968, 2969, 3014, 3019, 3020, 3021, 3022, 3031, 3032, 3033, 3092, 3093, 3094, 3096, 3097, 3099, 3102, 3104, 3113, 3114, 3115, 3124, 3137, 3138, 3142, 3143, 3879, 3881, 3896, 3898, 3901, 3902, 3915, 3917, 4262, 4264, 4275, 4276, 4291, 4292, 4293, 4294, 4295, 4296, 4297, 4298, 4299, 4300, 4301, 4302, 4303, 4304, 4305, 4307, 4312, 4313, 4314, 4315, 4316, 4317, 4318, 4320, 4328, 4330, 4333, 4334, 4530, 4531, 4532, 4533, 4535, 4538, 4539, 4540, 4541, 4547, 4548, 4549, 4550, 4551, 4552, 4555, 4556, 4560, 4561, 4564, 4565, 4570, 4571, 4572, 4573, 4822, 4823, 4824, 4825, 4826, 4827, 4829, 4831, 4833,
    4834, 4835, 4836, 4852, 4853, 4855, 4856, 4862, 4863, 4864, 4866, 4867, 4868, 4872, 4874, 4875, 4878
    \\[1.5mm]
    \multicolumn{1}{c}{\normalsize $14$ crossings} \\[1.5mm]
    432, 2459, 2651, 2757, 2762, 2904, 2936, 3115, 3147, 3228, 3233, 3238, 3249, 3264, 3279, 3291, 3292, 3299, 3305, 3311, 3741, 3755, 3757, 3779, 3780, 3911, 4515, 4517, 4518, 4541, 4542, 4708, 5249, 5250, 5251, 5286, 5287, 5392, 5609, 5689, 5763, 5834, 5889, 6048, 6127, 6166, 6211, 6432, 6436, 6438, 6440, 6477, 6478, 6479, 6485, 6488, 6489, 6493, 6503, 6511, 6521, 6522, 6528, 6533, 6535, 6537, 6543, 6551, 6555, 6557, 6561, 6564, 6565, 6574, 6579, 6582, 6689, 6690, 6691, 6708, 6711, 6712, 6721, 6728, 6733, 6734, 6739, 6740, 6751, 6753, 6755, 6756, 6759, 6778, 6781, 6782, 6788, 6790, 6793, 6794, 6802, 6803, 6806, 6900, 6916, 6944, 6945, 6948, 6949, 6951, 6955, 6958, 6959, 6965, 6966, 6967, 6973, 6974, 6975, 6977, 6979, 6985, 6997, 7020, 7027, 7032, 7054, 7061, 7063, 7077, 7091, 7094, 7107, 7109, 7113, 7120, 7122, 7124, 7126, 7322, 7323, 7324, 7371, 7373, 7375, 7376, 7378, 7382, 7385, 7386, 7388, 7392, 7393, 7394, 7435, 7518, 7551, 7552, 7554, 7555, 7556, 7557, 7558, 7563, 7564, 7565,
    7586, 7589, 7608, 7636, 7638, 7640, 7657, 7662, 7664, 7666, 7669, 7681, 7683, 7685, 7686, 7689, 7715, 7716, 7738, 8185, 8186, 8187, 8188, 8189, 8199, 8200, 8201, 8202, 8203, 8230, 8240, 8283, 8298, 8312, 8313, 8319, 8334, 8335, 8870, 8871, 8895, 8896, 8916, 8917, 8932, 8933, 9194, 9250, 9260, 9315, 9328, 9366, 9388, 9401, 9585, 9594, 9599, 10023, 10024, 10025, 10026, 10027, 10028, 10029, 10030, 10031, 10032, 10033, 10034, 10036, 10037, 10038, 10039, 10047, 10048, 10049, 10050, 10051, 10052, 10054, 10056, 10062, 10063, 10065, 10070, 10074, 10075, 10076, 10078, 10080, 10081, 10087, 10088, 10090, 10110, 10111, 10112, 10117, 10119, 10122, 10123, 10124, 10125, 10126, 10129, 10130, 10132, 10133, 10134, 10135, 10136, 10137, 10138, 10139, 10140, 10148, 10152, 10162, 10165, 10167, 10168, 10169, 10228, 10229, 10234, 10235, 10236, 10237, 10238, 10239, 10240, 10241, 10242, 10243, 10244, 10245, 10246, 10247, 10248, 10250, 10251, 10252, 10253, 10254, 10255, 10261, 10262, 10267, 10269, 10272,
    10273, 10274, 10325, 10327, 10333, 10335, 10336, 10338, 10339, 10341, 10342, 10374, 10376, 10380, 10382, 10383, 10384, 10385, 10511, 10512, 10513, 10514, 10515, 10516, 10518, 10519, 10520, 10521, 10522, 10523, 10524, 10525, 10526, 10527, 10528, 10529, 10539, 10540, 10550, 10551, 10552, 10553, 10554, 10555, 10561, 10562, 10632, 10633, 10634, 10635, 10641, 10643, 10644, 10645, 10646, 10647, 10648, 10653, 10654, 10655, 10656, 10694, 10695, 10698, 10699, 10700, 10701, 10704, 10992, 10993, 10994, 10995, 11009, 11010, 11011, 11012, 11013, 11014, 11015, 11017, 11018, 11019, 11080, 11081, 11085, 11086, 11087, 11088, 11104, 11105, 11424, 11425, 11447, 11448, 11452, 11453, 11462, 11463, 11479, 11480, 11643, 11654, 11655, 11685, 11689, 11690, 12180, 12181, 12182, 12184, 12185, 12186, 12187, 12188, 12189, 12191, 12192, 12194, 12195, 12196, 12197, 12198, 12199, 12201, 12202, 12204, 12205, 12206, 12207, 12208, 12209, 12210, 12211, 12212, 12216, 12219, 12222, 12231, 12233, 12234, 12235, 12236,
    12237, 12239, 12241, 12245, 12246, 12247, 12263, 12264, 12265, 12266, 12269, 12270, 12273, 12274, 12299, 12300, 12301, 12308, 12309, 12310, 12371, 12372, 12373, 12374, 12375, 12376, 12378, 12379, 12380, 12381, 12382, 12384, 12385, 12390, 12391, 12401, 12402, 12403, 12441, 12443, 12445, 12477, 12478, 12479, 12607, 12608, 12609, 12610, 12617, 12618, 12619, 12652, 12653, 12656, 12707, 12708, 12718, 12719, 12740, 12741, 14965, 14966, 15037, 15038, 15047, 15048, 15103, 15118, 16447, 16448, 16470, 16480, 16495, 16496, 16497, 16498, 16501, 16502, 16503, 16504, 16506, 16507, 16508, 16509, 16510, 16511, 16512, 16513, 16514, 16517, 16523, 16524, 16526, 16527, 16528, 16529, 16544, 16545, 16546, 16552, 16553, 16554, 16566, 16567, 16570, 16571, 16572, 16617, 16618, 16621, 16623, 16624, 16626, 16627, 16628, 16658, 16663, 17388, 17389, 17390, 17391, 17392, 17393, 17396, 17397, 17400, 17401, 17404, 17406, 17411, 17412, 17413, 17414, 17415, 17417, 17418, 17419, 17420, 17421, 17422, 17423, 17425,
    17426, 17427, 17428, 17429, 17431, 17432, 17434, 17435, 17436, 17438, 17440, 17441, 17443, 17445, 17449, 17450, 17451, 17452, 17453, 17467, 17468, 17469, 17470, 17471, 17472, 17473, 17474, 17475, 17485, 17486, 17487, 17509, 17510, 17517, 17518, 17701, 17702, 17703, 17705, 17706, 17708, 17711, 17713, 17715, 17716, 17717, 17719, 17723, 17724, 17729, 17730, 19298, 19299, 19300, 19301, 19302, 19303, 19304, 19305, 19306, 19307, 19308, 19309, 19312, 19313, 19314, 19315, 19316, 19320, 19321, 19322, 19323, 19326, 19327, 19328, 19329, 19364, 19365, 19366, 19367, 19369, 19370, 19371, 19373, 19374, 19420, 19421, 19422, 19424, 19425, 19426, 19428, 19429
  \end{tabular}
  \caption{%
    The $2$-bridge alternating prime knots $c^A_i$, for
    $c=8,\dots,14$. To these add all alternating prime knots of between
    $3$ and $7$ crossings.
  }%
  \label{table:Thetwobridgeknotsof3to14crossings}
\end{table}

\begin{table}[htbp]
  \tiny
  \hspace{-50pt}
  \begin{tabular}{p{150mm}}
    1621, 10492, 10680, 11059, 11100, 11532, 11555, 12252, 12256, 12677, 12678, 12815, 12816, 12840, 12864, 12883, 12913, 12933, 12938, 12980, 13028, 13041, 14448, 14449, 14769, 14817, 14822, 14827, 17009, 17010, 17427, 17484, 17497, 17499, 19699, 19700, 20014, 20089, 20091, 20092, 21965, 22016, 22017, 22018, 22559, 22776, 22824, 23043, 23105, 23243, 23292, 23432, 23467, 24846, 24852, 24857, 24876, 24978, 24989, 24991, 25009, 25100, 25111, 25112, 25120, 25135, 25137, 25141, 25154, 25164, 25173, 25174, 25189, 25190, 25194, 25212, 25215, 25228, 25235, 25241, 25253, 25255, 25262, 25291, 25302, 25317, 25690, 25721, 25723, 25731, 25833, 25838, 25839, 25846, 25872, 25873, 25874, 25883, 25884, 25901, 25919, 25920, 25926, 25927, 25938, 25946, 25956, 25958, 25961, 25962, 25969, 25986, 25987, 26013, 26086, 26103, 26120, 26554, 26561, 26562, 26566, 26584, 26585, 26586, 26591, 26592, 26604, 26605, 26631, 26632, 26634, 26635, 26645, 26646, 26647, 26648, 26654, 26658, 26659, 26663, 26675, 26676,
    26705, 26737, 26749, 26761, 26918, 26926, 26938, 26956, 26969, 26974, 26989, 27024, 27028, 27046, 27050, 27062, 27081, 27085, 27090, 27106, 27107, 27109, 27117, 27125, 27134, 27159, 27179, 27296, 27389, 27403, 27422, 27430, 27435, 28162, 28170, 28171, 28173, 28177, 28178, 28179, 28191, 28192, 28193, 28282, 28303, 28308, 28310, 28312, 28318, 28785, 28791, 28792, 28793, 28798, 28799, 28800, 28813, 28814, 28815, 28898, 28919, 28921, 28923, 28924, 28927, 29158, 29159, 29229, 29237, 29240, 29247, 29299, 29305, 29308, 29312, 29378, 29379, 29392, 29393, 29396, 29397, 29399, 29455, 29464, 29489, 29500, 29534, 29548, 29600, 29605, 29660, 29662, 29720, 29722, 31382, 31383, 31384, 31385, 31386, 31387, 31388, 31408, 31409, 31410, 31500, 31524, 31526, 31528, 31529, 31531, 31756, 31757, 31801, 31802, 31867, 31869, 31870, 31875, 31909, 31910, 31927, 31928, 31930, 31931, 31932, 31968, 31981, 32017, 32045, 32126, 32139, 32141, 32142, 32196, 32199, 32245, 32247, 33981, 33982, 34029, 34030, 34085,
    34086, 34127, 34128, 34153, 34154, 34155, 34156, 34157, 34206, 34282, 34333, 34344, 34414, 34428, 34458, 34459, 36024, 36092, 36115, 36116, 36170, 36190, 36191, 36233, 36246, 36247, 37138, 37322, 37343, 37421, 37457, 37489, 37559, 37585, 37620, 37659, 39530, 39531, 39532, 39533, 39538, 39548, 39549, 39551, 39553, 39559, 39570, 39587, 39604, 39606, 39609, 39611, 39614, 39615, 39617, 39621, 39622, 39623, 39629, 39631, 39637, 39638, 39639, 39641, 39647, 39651, 39652, 39654, 39655, 39658, 39659, 39667, 39668, 39670, 39674, 39681, 39682, 39684, 39686, 39689, 39690, 39694, 39696, 39699, 39700, 39702, 39706, 39735, 39742, 39849, 39850, 39874, 39876, 39879, 39880, 39881, 39882, 39887, 39888, 39889, 39893, 39894, 39899, 39900, 39901, 39902, 39903, 39904, 39905, 39908, 39909, 39913, 39914, 39916, 39917, 39922, 39923, 39924, 39925, 39926, 39927, 39928, 39931, 39932, 39936, 39937, 39939, 39940, 39944, 39945, 39946, 39949, 39961, 39968, 39969, 39975, 39995, 39999, 40018, 40026, 40027, 40041,
    40046, 40049, 40063, 40067, 40071, 40072, 40075, 40077, 40081, 40092, 40093, 40095, 40096, 40099, 40102, 40105, 40106, 40120, 40128, 40144, 40162, 40165, 40252, 40256, 40520, 40521, 40537, 40538, 40539, 40540, 40541, 40542, 40543, 40547, 40558, 40559, 40560, 40561, 40562, 40563, 40564, 40565, 40578, 40579, 40581, 40582, 40585, 40586, 40587, 40588, 40589, 40590, 40591, 40593, 40594, 40595, 40596, 40600, 40608, 40609, 40610, 40615, 40621, 40635, 40666, 40667, 40677, 40694, 40696, 40698, 40699, 40701, 40704, 40705, 40710, 40711, 40712, 40713, 40715, 40717, 40719, 40720, 40724, 40725, 40727, 40730, 40731, 40820, 41088, 41089, 41093, 41094, 41099, 41105, 41106, 41110, 41111, 41113, 41114, 41117, 41122, 41123, 41125, 41126, 41133, 41134, 41135, 41136, 41137, 41140, 41141, 41143, 41147, 41335, 41345, 41351, 41353, 41356, 41399, 41413, 41417, 41565, 41568, 41570, 42067, 42068, 42069, 42070, 42071, 42072, 42075, 42076, 42077, 42078, 42079, 42080, 42085, 42086, 42087, 42088, 42123, 42128,
    42135, 42136, 42197, 42198, 42220, 42221, 42223, 42224, 42225, 42226, 42229, 42230, 42231, 42250, 42251, 42253, 42256, 42257, 42258, 42259, 42260, 42263, 42264, 42293, 42336, 42339, 42579, 42580, 42581, 42600, 42601, 42604, 42605, 42606, 42607, 42608, 42616, 42617, 42618, 42622, 42623, 42624, 42625, 42626, 42627, 42628, 42629, 42637, 42639, 42640, 42645, 42646, 42684, 42685, 42829, 42843, 42850, 42855, 42866, 42867, 42888, 42890, 42900, 42901, 43059, 43061, 43063, 43064, 44042, 44043, 44044, 44063, 44064, 44065, 44083, 44084, 44085, 44101, 44102, 44103, 44104, 44105, 44106, 44108, 44109, 44124, 44125, 44177, 44178, 44365, 44371, 44374, 44384, 44385, 44416, 44417, 44423, 44424, 44425, 44426, 44437, 44439, 44611, 44612, 44623, 44624, 45915, 45916, 45944, 45945, 45955, 45956, 45957, 45958, 45988, 45989, 45996, 45997, 45998, 45999, 46260, 46261, 46269, 46270, 47303, 47368, 47378, 47382, 47383, 47414, 47426, 47429, 47430, 47442, 47445, 47446, 47928, 47938, 47941, 47950, 47954, 47965,
    50203, 50204, 50205, 50206, 50207, 50208, 50209, 50210, 50211, 50214, 50215, 50216, 50217, 50218, 50219, 50225, 50226, 50227, 50228, 50229, 50230, 50231, 50232, 50234, 50235, 50236, 50237, 50238, 50240, 50241, 50244, 50245, 50246, 50253, 50254, 50255, 50257, 50258, 50259, 50260, 50261, 50262, 50263, 50264, 50267, 50268, 50269, 50270, 50271, 50272, 50274, 50275, 50276, 50279, 50280, 50285, 50286, 50289, 50295, 50296, 50297, 50298, 50299, 50309, 50310, 50311, 50312, 50313, 50315, 50316, 50318, 50320, 50321, 50322, 50323, 50324, 50325, 50326, 50331, 50355, 50379, 50386, 50391, 50392, 50393, 50394, 50395, 50396, 50397, 50399, 50401, 50402, 50403, 50404, 50406, 50407, 50408, 50409, 50410, 50411, 50412, 50413, 50414, 50450, 50451, 50452, 50453, 50454, 50455, 50456, 50457, 50460, 50486, 50522, 50594, 50595, 50596, 50597, 50599, 50600, 50601, 50602, 50603, 50604, 50606, 50607, 50608, 50609, 50610, 50611, 50612, 50614, 50615, 50616, 50617, 50618, 50619, 50620, 50643, 50644, 50646, 50647,
    50649, 50650, 50651, 50652, 50767, 50771, 50772, 50773, 50774, 50775, 50776, 50777, 50779, 50780, 50794, 50795, 50873, 50875, 50876, 50877, 50902, 50909, 51192, 51193, 51194, 51195, 51196, 51197, 51198, 51199, 51200, 51201, 51205, 51206, 51207, 51209, 51210, 51211, 51212, 51213, 51215, 51216, 51217, 51218, 51219, 51220, 51222, 51223, 51225, 51226, 51227, 51228, 51229, 51230, 51231, 51232, 51233, 51246, 51249, 51274, 51275, 51276, 51277, 51279, 51280, 51282, 51283, 51284, 51291, 51292, 51297, 51299, 51327, 51328, 51435, 51436, 51437, 51438, 51439, 51441, 51442, 51443, 51447, 51448, 51449, 51457, 51458, 51461, 51462, 51464, 51553, 51554, 51559, 51563, 51595, 51596, 51693, 51694, 52180, 52181, 52182, 52193, 52194, 52195, 52215, 52216, 52217, 52218, 52219, 52220, 52221, 52222, 52223, 52225, 52226, 52240, 52241, 52278, 52279, 52392, 52396, 52399, 52420, 52421, 52442, 52444, 52566, 52567, 53093, 53094, 53112, 53113, 53115, 53116, 53117, 53118, 53119, 53120, 53217, 53218, 53414, 53430,
    53431, 53432, 53433, 53434, 53509, 53511, 53512, 53513, 54584, 54585, 54586, 54588, 54589, 54591, 54592, 54593, 54597, 54602, 54604, 54606, 54613, 54614, 54625, 54627, 54670, 54671, 54672, 54674, 54675, 54689, 54730, 54752, 54824, 54825, 54826, 54850, 54880, 54881, 54893, 54894, 65122, 65131, 65514, 65519, 65527, 65528, 65677, 65682, 70803, 70808, 71021, 71022, 71315, 71317, 71318, 71326, 71343, 71344, 71345, 71349, 71350, 71351, 71355, 71356, 71357, 71359, 71360, 71363, 71364, 71374, 71388, 71390, 71394, 71395, 71398, 71399, 71400, 71402, 71424, 71426, 71458, 71459, 71460, 71465, 71466, 71467, 71495, 71497, 71501, 71502, 71506, 71507, 71509, 71510, 71511, 71512, 71525, 71526, 71603, 71604, 71608, 71609, 71615, 71616, 71617, 71644, 71645, 71739, 71741, 71751, 71753, 71764, 71768, 71850, 71852, 71853, 71854, 71859, 71862, 75697, 75699, 75700, 75705, 75720, 75721, 75722, 75723, 75724, 75727, 75728, 75730, 75736, 75746, 75748, 75779, 75780, 75786, 75787, 75826, 75827, 75828, 75841,
    75985, 75986, 75987, 75988, 75989, 75990, 75991, 75992, 75993, 75994, 75995, 75996, 75997, 75998, 75999, 76000, 76001, 76002, 76004, 76005, 76006, 76007, 76008, 76009, 76012, 76016, 76017, 76018, 76019, 76020, 76021, 76022, 76023, 76024, 76025, 76029, 76040, 76041, 76042, 76043, 76044, 76045, 76046, 76047, 76048, 76049, 76065, 76066, 76067, 76068, 76072, 76077, 76113, 76114, 76115, 76116, 76117, 76118, 76119, 76120, 76121, 76122, 76123, 76125, 76126, 76131, 76141, 76142, 76144, 76146, 76160, 76161, 76221, 76222, 76223, 76224, 76231, 76232, 76236, 76237, 76238, 76239, 76240, 76255, 76286, 76288, 76289, 76290, 76295, 76296, 76297, 76298, 78550, 78551, 78552, 78553, 78554, 78555, 78556, 78557, 78558, 78562, 78563, 78564, 78566, 78572, 78573, 78581, 78586, 78597, 78598, 78599, 78600, 78601, 78603, 78604, 78606, 78609, 78611, 78616, 78617, 78644, 78645, 78646, 78653, 78720, 78721, 78722, 78724, 78725, 78726, 78727, 78728, 78730, 78731, 78734, 78735, 78743, 78745, 78754, 78755, 78796,
    78797, 78798, 78802, 78803, 78804, 78820, 78822, 78853, 78854, 78855, 78856, 78864, 78865, 78879, 78880, 84276, 84277, 84278, 84279, 84280, 84283, 84284, 84285, 84286, 84287, 84288, 84289, 84290, 84291, 84292, 84296, 84298, 84300, 84301, 84302, 84303, 84718, 84719, 84769, 84771, 84772, 84773, 84774, 84775, 84776, 84777, 84778, 84779, 84780, 84784, 84785, 84786, 84789, 84790, 84793, 84794, 84797, 84798, 84805, 84806, 84807, 84808, 84809, 84810, 84956, 84957, 84958, 84959, 84961, 84964, 84965, 84968, 84969, 85170, 85171, 85172, 85173, 85174, 85175, 85178, 85179, 85180, 85182, 85184, 85186, 85190, 85192, 85193, 85194, 85195, 85197, 85198, 85204, 85205, 85207, 85208, 85226, 85228, 85234, 85240, 85241, 85242, 85245, 85246, 85247, 85252, 85258, 85259, 85263
  \end{tabular}
  \caption{%
    Identified $2$-bridge alternating prime knots $15^A_i$.
  }%
\end{table}

\begin{table}[htbp]
  \tiny
  \hspace{-50pt}
  \begin{tabular}{p{150mm}}
    7016, 46278, 48134, 49567, 49587, 51133, 51374, 53593, 53959, 55922, 56092, 56810, 56830, 56850, 56908, 57069, 57228, 57305, 57385, 57386, 57423, 57454, 57479, 57501, 63516, 63643, 63648, 63789, 63790, 64888, 73048, 73065, 73067, 73302, 73303, 74931, 83426, 83428, 83429, 83663, 83664, 85085, 92484, 92485, 92486, 92726, 92727, 93380, 95964, 96750, 97688, 98379, 99037, 99608, 101329, 101978, 102657, 102957, 103322, 107661, 107694, 107699, 107704, 108045, 108081, 108083, 108111, 108473, 108474, 108475, 108566, 108605, 108606, 108635, 108688, 108759, 108791, 108867, 108868, 108892, 108913, 108918, 108923, 108942, 108993, 109017, 109029, 109051, 109063, 109096, 109107, 109108, 109148, 109184, 109212, 109219, 110462, 110472, 110476, 110537, 110895, 110896, 110897, 111001, 111012, 111013, 111055, 111075, 111130, 111168, 111169, 111239, 111240, 111309, 111318, 111328, 111330, 111349, 111458, 111469, 111471, 111505, 111511, 111548, 111559, 111609, 111610, 111667, 111703, 111706, 111728,
    113388, 113389, 113390, 113513, 113516, 113517, 113576, 113585, 113637, 113638, 113658, 113659, 113701, 113702, 113780, 113784, 113796, 113797, 113828, 113908, 113919, 113920, 113954, 113957, 113958, 113991, 113993, 114025, 114026, 114055, 114056, 114064, 114099, 115013, 115194, 115460, 115672, 115673, 115682, 115683, 115688, 115702, 115711, 115712, 115739, 115742, 115743, 115772, 115773, 115774, 115801, 115802, 115803, 115814, 115822, 115829, 115856, 115959, 116077, 116208, 116240, 116257, 116394, 116446, 116458, 116546, 116691, 116704, 116802, 116818, 116877, 116908, 116975, 116987, 117005, 117045, 117050, 117064, 117087, 117114, 117129, 117139, 117148, 117154, 120115, 120116, 120117, 120433, 120438, 120467, 120469, 120475, 120499, 120510, 120511, 120519, 120534, 120536, 120540, 120806, 122007, 122008, 122009, 122309, 122311, 122318, 122319, 122322, 122357, 122362, 122363, 122370, 122396, 122397, 122398, 122654, 123548, 123713, 123976, 123977, 124002, 124003, 124005, 124023,
    124030, 124031, 124035, 124053, 124054, 124055, 124218, 124257, 124445, 124632, 124834, 124839, 124844, 124980, 125016, 125027, 125035, 125047, 125180, 125193, 125215, 125220, 125235, 125393, 125409, 125421, 125425, 125437, 125681, 125682, 125881, 125995, 126089, 131793, 131794, 131795, 132170, 132172, 132174, 132175, 132177, 132221, 132229, 132230, 132232, 132236, 132237, 132238, 132602, 133607, 133836, 133838, 133840, 133841, 133842, 133886, 133891, 133892, 133893, 133898, 133899, 133900, 134057, 134103, 134246, 134554, 134566, 134568, 134729, 134731, 134733, 134754, 134755, 134891, 134896, 134933, 134941, 134944, 135041, 135043, 135075, 135080, 135083, 135341, 135342, 135506, 135638, 136003, 143292, 143293, 143294, 143295, 143296, 143375, 143376, 143377, 143378, 143379, 143380, 143381, 143534, 143594, 144076, 144078, 144079, 144217, 144219, 144220, 144243, 144244, 144431, 144432, 144433, 144460, 144461, 144542, 144543, 144583, 144585, 144586, 144700, 144701, 144840, 145294,
    145564, 151696, 151697, 151698, 151919, 151920, 151921, 151957, 151958, 152142, 152143, 152144, 152183, 152184, 152342, 152343, 152344, 152376, 152377, 152487, 152488, 153052, 153200, 153437, 156081, 156678, 156737, 157220, 157308, 157403, 157737, 157838, 157955, 158069, 158499, 158663, 158825, 158912, 161248, 161606, 161676, 162086, 162170, 162407, 162492, 162560, 162616, 169837, 169839, 169843, 169845, 169854, 169856, 169862, 169865, 169966, 169977, 169978, 169979, 169981, 170000, 170001, 170002, 170004, 170009, 170010, 170011, 170012, 170029, 170031, 170034, 170035, 170047, 170048, 170051, 170052, 170056, 170066, 170086, 170090, 170098, 170117, 170119, 170122, 170123, 170133, 170134, 170136, 170152, 170153, 170155, 170159, 170160, 170161, 170167, 170172, 170190, 170193, 170195, 170201, 170215, 170245, 170249, 170256, 170260, 170264, 170278, 170284, 170285, 170320, 170339, 170342, 170352, 170740, 170746, 170747, 170748, 170750, 170763, 170764, 170765, 170771, 170772, 170773,
    170774, 170775, 170787, 170788, 170792, 170795, 170796, 170804, 170805, 170807, 170808, 170813, 170820, 170827, 170828, 170829, 170830, 170835, 170836, 170845, 170849, 170850, 170855, 170856, 170858, 170868, 170869, 170870, 170875, 170876, 170877, 170887, 170889, 170902, 170908, 170909, 170912, 170937, 170960, 170962, 170965, 170966, 170972, 170975, 170977, 170983, 170987, 170988, 171003, 171004, 171021, 171033, 171110, 171303, 171310, 171312, 171319, 171322, 171323, 171350, 171355, 171359, 171362, 171368, 171372, 171376, 171391, 171395, 171396, 171431, 171438, 171441, 171454, 171457, 171466, 171467, 171476, 171477, 171480, 171481, 171483, 171496, 171497, 171520, 171523, 171524, 171536, 171540, 171541, 171542, 171551, 171557, 171558, 171559, 171582, 171593, 171598, 171599, 171617, 171626, 171635, 171684, 171694, 171703, 171707, 171769, 171781, 171787, 171797, 171799, 171832, 171844, 171913, 171930, 171937, 171971, 171999, 172004, 172029, 172038, 172057, 172060, 172064, 172069,
    172071, 172078, 172084, 173180, 173181, 173182, 173204, 173205, 173206, 173207, 173208, 173209, 173210, 173211, 173248, 173249, 173250, 173258, 173259, 173261, 173262, 173285, 173289, 173299, 173301, 173302, 173303, 173311, 173312, 173313, 173314, 173315, 173316, 173317, 173318, 173336, 173337, 173338, 173339, 173340, 173341, 173361, 173363, 173370, 173373, 173374, 173376, 173386, 173412, 173414, 173418, 173419, 173421, 173432, 173433, 173436, 173443, 173444, 173459, 173460, 173462, 173472, 173718, 173719, 173726, 173727, 173733, 173736, 173737, 173754, 173755, 173770, 173771, 173776, 173785, 173786, 173787, 173788, 173805, 173808, 173809, 173810, 173815, 173822, 173823, 173837, 173838, 173843, 173844, 173846, 173847, 173848, 173852, 173853, 173873, 173877, 173878, 173879, 173886, 173887, 173889, 173894, 173895, 173896, 173909, 173910, 173911, 173915, 173916, 173929, 173936, 173947, 173965, 173975, 173982, 173984, 173987, 173988, 173990, 174010, 174024, 174025, 174041, 174047,
    174071, 174088, 174095, 174127, 174142, 174150, 174177, 174180, 174187, 174189, 174191, 174193, 175407, 175408, 175415, 175416, 175429, 175432, 175433, 175453, 175461, 175462, 175465, 175466, 175471, 175472, 175473, 175504, 175525, 175526, 175528, 175529, 175534, 175535, 175536, 175540, 175541, 175553, 175557, 175558, 175559, 175568, 175569, 175570, 175571, 175584, 175585, 175586, 175588, 175589, 175592, 175593, 175597, 175598, 175613, 175614, 175615, 175616, 175644, 175645, 175674, 175684, 175685, 175698, 175753, 175793, 175796, 175804, 175806, 175811, 175812, 175816, 176488, 176509, 176577, 176588, 176605, 176620, 176665, 176673, 176689, 176695, 176702, 176743, 176755, 176758, 176780, 176781, 176812, 176814, 176836, 176840, 176848, 176849, 176861, 176862, 176866, 176869, 176873, 177566, 177572, 177592, 177599, 177602, 177619, 177621, 177630, 177633, 177650, 177652, 177655, 179695, 179696, 179697, 179698, 179699, 179700, 179701, 179702, 179703, 179704, 179831, 179833, 179839,
    179842, 179843, 179845, 179850, 179854, 179856, 179862, 179863, 179865, 179869, 179870, 179871, 179877, 179884, 179885, 179989, 180175, 180176, 180177, 180238, 180239, 180240, 180241, 180242, 180243, 180255, 180256, 180321, 180322, 180325, 180326, 180328, 180329, 180330, 180333, 180334, 180337, 180338, 180339, 180344, 180345, 180346, 180350, 180355, 180356, 180357, 180358, 180359, 180409, 180412, 180416, 180428, 180467, 180469, 180472, 180481, 180503, 180531, 180535, 180537, 180539, 180545, 180587, 180594, 180597, 180602, 180605, 180646, 180651, 180656, 180660, 180661, 180664, 180666, 180670, 180728, 180729, 180783, 180786, 181784, 181785, 181840, 181841, 181843, 181844, 181847, 181848, 181849, 181850, 181851, 181855, 181856, 181857, 181858, 181859, 181860, 181864, 181875, 181876, 181877, 181878, 181879, 181917, 181918, 181919, 181923, 181963, 181964, 182005, 182070, 182072, 182078, 182080, 182082, 182085, 182105, 182106, 182676, 182754, 182783, 182813, 182908, 182949, 182952,
    182953, 182955, 182983, 182987, 182990, 182991, 182995, 182996, 183001, 183018, 183020, 183024, 183025, 183027, 183028, 183031, 183083, 183084, 183085, 183118, 183130, 183438, 183775, 183777, 183793, 183797, 183802, 183816, 183818, 183825, 183848, 183860, 183870, 188126, 188127, 188129, 188130, 188131, 188132, 188133, 188135, 188136, 188163, 188164, 188165, 188166, 188167, 188168, 188171, 188172, 188173, 188174, 188175, 188176, 188239, 188241, 188242, 188265, 188377, 188415, 188443, 188445, 188456, 188458, 188459, 188460, 188502, 188503, 189521, 189523, 189524, 189525, 189571, 189572, 189573, 189574, 189579, 189580, 189581, 189598, 189600, 189609, 189610, 189613, 189614, 189615, 189646, 189647, 189648, 189675, 190218, 190220, 190244, 190245, 190248, 190267, 190268, 190276, 190287, 190304, 190403, 190404, 196138, 196139, 196140, 196141, 196207, 196208, 196209, 196210, 196219, 196220, 196221, 196288, 196289, 196290, 196291, 196302, 196303, 196304, 196336, 196337, 196338, 197100,
    197101, 197127, 197128, 197134, 197149, 197150, 197153, 197160, 197251, 197252, 197330, 197331, 202220, 202221, 202396, 202397, 202417, 202418, 202543, 202544, 202568, 202569, 202590, 202591, 202700, 202701, 202757, 202758, 202827, 202828, 206228, 206324, 206351, 206440, 206473, 206490, 207872, 207898, 207904, 207929, 207941, 207950, 219644, 219645, 219646, 219647, 219648, 219649, 219650, 219651, 219652, 219653, 219654, 219658, 219659, 219660, 219661, 219662, 219663, 219664, 219667, 219668, 219669, 219670, 219671, 219672, 219673, 219674, 219675, 219676, 219677, 219679, 219681, 219682, 219685, 219686, 219687, 219689, 219690, 219691, 219692, 219700, 219701, 219710, 219711, 219712, 219713, 219714, 219715, 219716, 219717, 219718, 219719, 219720, 219722, 219723, 219724, 219725, 219726, 219727, 219728, 219729, 219730, 219731, 219733, 219735, 219741, 219751, 219752, 219754, 219759, 219763, 219771, 219779, 219780, 219781, 219782, 219784, 219788, 219792, 219793, 219808, 219809, 219811,
    219816, 219879, 219880, 219881, 219882, 219884, 219885,
    219886, 219888, 219889, 219895, 219900, 219901, 219902, 219903, 219904, 219905, 219906, 219908, 219911, 219912, 219919, 219920, 219921, 219922, 219923, 219926, 219927, 219932, 219933, 219939, 219940, 219941, 219944, 219945, 219946, 219947, 219953, 219954, 219955, 219957, 219958, 219980, 219984, 220023, 220032, 220038, 220039, 220040, 220045, 220228, 220229, 220230, 220231, 220232, 220233, 220234, 220235, 220236, 220237, 220240, 220241, 220242, 220243, 220244, 220245, 220246, 220250, 220251, 220252, 220257, 220258, 220259, 220260, 220261, 220262, 220263, 220264, 220268, 220269, 220270, 220271, 220272, 220273, 220274, 220275, 220276, 220277, 220278, 220279, 220281, 220282, 220283, 220284, 220285, 220286, 220287, 220288, 220289, 220291, 220292, 220295, 220299, 220300, 220302, 220306, 220307, 220308, 220325, 220330, 220336, 220337, 220338, 220474, 220475, 220477, 220482, 220486, 220495, 220500, 220502, 220507, 220510, 220513, 220518, 220525, 220539, 220543, 220545, 220552, 220553,
    220554, 220556, 220557, 220681, 220684, 220689, 220691,
  \end{tabular}
  \caption{%
    Identified $2$-bridge alternating prime knots $16^A_i$, part 1/2.
  }%
\end{table}

\begin{table}[htbp]
  \tiny
  \hspace{-50pt}
  \begin{tabular}{p{150mm}}
    220694, 220695, 220696, 220697, 221072, 221073, 221074, 221075, 221076, 221077, 221078, 221079, 221080, 221081, 221085, 221086, 221087, 221094, 221095, 221096, 221098, 221099, 221100, 221101, 221102, 221103, 221105, 221106, 221107, 221108, 221109, 221110, 221111, 221112, 221113, 221114, 221115, 221117, 221118, 221119, 221120, 221121, 221123, 221124, 221126, 221127, 221128, 221129, 221130, 221131, 221132, 221133, 221134, 221135, 221136, 221137, 221138, 221139, 221148, 221156, 221157, 221158, 221166, 221170, 221174, 221175, 221176, 221271, 221272, 221273, 221274, 221275, 221278, 221279, 221280, 221281, 221282, 221283, 221285, 221287, 221290, 221291, 221294, 221295, 221303, 221304, 221308, 221310, 221312, 221314, 221316, 221391, 221392, 221396, 221397, 221398, 221823, 221824, 221829, 221830, 221835, 221836, 221839, 221840, 221842, 221843, 221847, 221848, 221865, 221866, 221869, 221870, 221871, 221872, 221876, 221877, 221878,
    221879, 221882, 221883, 221884, 221885, 221984, 221985, 222346, 222347, 222348, 222370, 222373, 222378, 222382, 222383, 222384, 222385, 222386, 222387, 222595, 222596, 222597, 222600, 222602, 222606, 222608, 222609, 222610, 222611, 223836, 223837, 223838, 223839, 223896, 223897, 223898, 223900, 223901, 223902, 223903, 223904, 223905, 223907, 223908, 223910, 223911, 223912, 223913, 223914, 223915, 223916, 223917, 223918, 223919, 223920, 223921, 223922, 223923, 223924, 223925, 223929, 223930, 223931, 223932, 223933, 223934, 223935, 223936, 223938, 223939, 223940, 223957, 223958, 223959, 223960, 223961, 223962, 223963, 223964, 223965, 223966, 223969, 223970, 223971, 223999, 224002, 224005, 224006, 224007, 224008, 224009, 224011, 224012, 224025, 224026, 224187, 224188, 224201, 224203, 224208, 224214, 224215, 224217, 224227, 224230, 224231, 224232, 224233, 224234, 224237, 224238, 224244, 224246, 224247, 224248, 224249, 224250, 224251, 224252, 224253, 224255, 224262, 224263, 224264,
    224265, 224266, 224267, 224359, 224362, 224366, 224367, 224368, 224369, 224370, 224371, 224373, 224376, 224915, 224916, 224926, 224927, 224938, 224939, 224940, 224941, 224942, 224943, 224944, 224945, 224959, 224960, 224963, 224964, 224965, 224966, 224967, 224968, 224970, 224971, 224972, 224973, 224974, 224975, 224982, 224983, 224984, 224986, 224987, 224988, 224989, 224997, 224998, 225024, 225025, 225067, 225068, 225074, 225075, 225077, 225078, 225081, 225471, 225473, 225497, 225506, 225507, 225513, 225514, 225515, 225516, 225517, 225518, 225519, 225521, 225524, 225966, 225968, 225973, 225981, 225982, 225985, 225986, 225987, 225988, 225989, 228502, 228503, 228504, 228505, 228506, 228507, 228509, 228510, 228545, 228546, 228547, 228548, 228549, 228550, 228552, 228553, 228554, 228555, 228556, 228557, 228558, 228559, 228560, 228561, 228562, 228563, 228568, 228569, 228570, 228572, 228573, 228592, 228593, 228616, 228617, 228706, 228707, 228718, 228719, 228720, 228721, 228722, 228723,
    228726, 228755, 228756, 229353, 229354, 229355, 229356, 229393, 229395, 229403, 229405, 229406, 229407, 229408, 229409, 229410, 229412, 229448, 229449, 229450, 229484, 229901, 229902, 229914, 229915, 229919, 229920, 229921, 229923, 229932, 229948, 233614, 233615, 233616, 233617, 233680, 233681, 233682, 233683, 233690, 233691, 233692, 233693, 233694, 233695, 233696, 233727, 233728, 233729, 233767, 233768, 233769, 234319, 234320, 234339, 234340, 234341, 234342, 234343, 234351, 234358, 234477, 234478, 237005, 237006, 237099, 237100, 237109, 237110, 237115, 237116, 237202, 237203, 237245, 237246, 237259, 237260, 238315, 238363, 238369, 238370, 238613, 238619, 238624, 238625, 249126, 249127, 249128, 249130, 249131, 249132, 249136, 249137, 249138, 249140, 249141, 249142, 249143, 249144, 249146, 249147, 249149, 249150, 249151, 249156, 249157, 249158, 249160, 249161, 249162, 249163, 249164, 249166, 249167, 249169, 249170, 249171, 249176, 249177, 249178, 249180, 249181, 249182, 249183,
    249184, 249185, 249186, 249187, 249190, 249191, 249192, 249193, 249194, 249195, 249197, 249198, 249204, 249209, 249214, 249219, 249220, 249235, 249236, 249237, 249239, 249240, 249241, 249242, 249243, 249244, 249245, 249248, 249253, 249257, 249270, 249273, 249276, 249310, 249313, 249315, 249316, 249317, 249318, 249320, 249321, 249322, 249325, 249326, 249329, 249330, 249362, 249364, 249365, 249366, 249367, 249384, 249387, 249388, 249416, 249417, 249418, 249482, 249483, 249484, 249485, 249486, 249487, 249488, 249492, 249493, 249494, 249498, 249499, 249500, 249547, 249548, 249549, 249551, 249553, 249657, 249658, 249681, 249682, 249683, 249684, 249685, 249686, 249781, 249782, 249783, 249793, 249794, 249795, 250406, 250407, 250408, 250410, 250411, 250412, 250413, 250414, 250415, 250417, 250418, 250420, 250421, 250422, 250423, 250424, 250425, 250427, 250428, 250430, 250431, 250432, 250445, 250446, 250447, 250448, 250449, 250450, 250474, 250477, 250480, 250555, 250557, 250558, 250563,
    250564, 250576, 250578, 250580, 250634, 250635, 250636, 250840, 250841, 250856, 250857, 250860, 250861, 250864, 250865, 251128, 251129, 251130, 251314, 251315, 251316, 252417, 252418, 252419, 252420, 252421, 252422, 252424, 252425, 252442, 252443, 252444, 252446, 252447, 252464, 252465, 252512, 252513, 252518, 252776, 252778, 252783, 252976, 252977, 252980, 253856, 253857, 253858, 253859, 253889, 253890, 253891, 254013, 254014, 254024, 254331, 254332, 254368, 254369, 254509, 254510, 291474, 291475, 292837, 292838, 292873, 292874, 293449, 293758, 312699, 312700, 313514, 313781, 314427, 314437, 314440, 314441, 314482, 314504, 314506, 314507, 314559, 314562, 314563, 314565, 314571, 314574, 314575, 314577, 314582, 314583, 314586, 314587, 314593, 314594, 314595, 314603, 314607, 314652, 314654, 314660, 314662, 314666, 314668, 314670, 314872, 314874, 314875, 314876, 314889, 314891, 314892, 314893, 314898, 314899, 314900, 314953, 314955, 314962, 314964, 314966, 314989, 314990, 315311,
    315312, 315313, 315314, 315335, 315336, 315337, 315338, 315343, 315344, 315345, 315437, 315438, 315446, 315447, 315450, 315466, 315467, 315958, 315959, 315968, 315969, 315983, 315984, 315993, 315994, 316009, 316010, 316020, 316021, 316032, 316033, 316443, 316459, 316468, 316620, 316624, 316627, 333177, 333183, 333186, 333187, 333201, 333209, 333211, 333212, 333239, 333240, 333243, 333244, 333250, 333251, 333252, 333254, 333256, 333281, 333283, 333285, 333356, 333357, 333358, 333379, 333384, 333385, 333474, 333475, 333476, 333501, 333507, 333508, 333582, 333583, 333696, 333777, 333779, 333783, 334022, 334023, 334024, 334025, 334027, 334028, 334029, 334030, 334031, 334032, 334033, 334034, 334035, 334036, 334037, 334045, 334046, 334048, 334049, 334051, 334053, 334055, 334056, 334057, 334058, 334059, 334062, 334063, 334065, 334067, 334068, 334069, 334070, 334073, 334074, 334075, 334076, 334077, 334078, 334079, 334080, 334081, 334082, 334083, 334097, 334099, 334101, 334103, 334104,
    334105, 334106, 334134, 334135, 334136, 334137, 334141, 334142, 334143, 334144, 334145, 334146, 334147, 334151, 334152, 334153, 334154, 334155, 334156, 334157, 334158, 334159, 334160, 334161, 334162, 334163, 334165, 334166, 334170, 334171, 334177, 334178, 334181, 334182, 334187, 334188, 334192, 334193, 334196, 334197, 334236, 334237, 334240, 334241, 334244, 334245, 334283, 334284, 334286, 334287, 334288, 334289, 334388, 334389, 334390, 334391, 334392, 334394, 334395, 334409, 334410, 334411, 334412, 334413, 334415, 334416, 334417, 334418, 334419, 334420, 334421, 334443, 334444, 334447, 334448, 334449, 334450, 334530, 334532, 334536, 334538, 334539, 334544, 334545, 334546, 334584, 334585, 334588, 334589, 334590, 334591, 334931, 334932, 334933, 334939, 334940, 334941, 334946, 334947, 334948, 334953, 334954, 334955, 334958, 334959, 334960, 335033, 335037, 335039, 335041, 335042, 335437, 335438, 335448, 335449, 335581, 335582, 335603, 335604, 344673, 344674, 344675, 344676, 344677,
    344680, 344681, 344682, 344683, 344684, 344688, 344689, 344690, 344697, 344698, 344699, 344700, 344701, 344702, 344704, 344705, 344706, 344707, 344711, 344712, 344727, 344728, 344738, 344739, 344756, 344757, 344758, 344759, 344760, 344761, 344771, 344772, 344773, 344774, 344780, 344781, 344792, 344793, 344794, 344798, 344824, 344825, 344829, 344830, 344853, 344854, 344865, 344866, 345045, 345046, 345047, 345048, 345049, 345050, 345051, 345052, 345054, 345055, 345057, 345058, 345059, 345060, 345061, 345062, 345063, 345064, 345065, 345066, 345067, 345069, 345070, 345072, 345073, 345074, 345075, 345076, 345077, 345079, 345080, 345082, 345083, 345084, 345085, 345086, 345087, 345088, 345089, 345090, 345095, 345099, 345102, 345109, 345111, 345112, 345114, 345115, 345116, 345117, 345118, 345120, 345121, 345123, 345125, 345129, 345130, 345131, 345132, 345133, 345152, 345153, 345156, 345157, 345158, 345159, 345162, 345163, 345166, 345167, 345203, 345204, 345205, 345206, 345207, 345223,
    345224, 345225, 345226, 345227, 345297, 345298, 345299, 345300, 345301, 345303, 345304, 345307, 345308, 345309, 345310, 345311, 345313, 345314, 345315, 345316, 345317, 345319, 345320, 345336, 345337, 345338, 345387, 345389, 345391, 345393, 345395, 345450, 345451, 345453, 345454, 345455, 345616, 345617, 345618, 345619, 345620, 345621, 345626, 345627, 345628, 345663, 345664, 345667, 345744, 345745, 345787, 345804, 350001, 350002, 350003, 350005, 350006, 350008, 350009, 350010, 350014, 350019, 350021, 350023, 350033, 350034, 350046, 350048, 350064, 350065, 350066, 350068, 350069, 350071, 350077, 350082, 350118, 350119, 350120, 350131, 350155, 350156, 350172, 350194, 372041, 372042, 372043, 372044, 372045, 372046, 372051, 372052, 372066, 372067, 372068, 372072, 372073, 372074, 372078, 372079, 372080, 372085, 372087, 372089, 372108, 372109, 372110, 372112, 372114, 372115, 372116, 374259, 374260, 374261, 374264, 374409, 374410, 374411, 374413, 374415, 374416, 374417, 374418, 374419,
    374421, 374422, 374423, 374424, 374425, 374427, 374428, 374429, 374430, 374431, 374433, 374434, 374435, 374436, 374437, 374439, 374440, 374441, 374442, 374443, 374445, 374446, 374449, 374450, 374451, 374452, 374453, 374457, 374459, 374461, 374463, 374465, 374471, 374472, 374473, 374474, 374475, 374476, 374477, 374478, 374493, 374494, 374495, 374496, 374497, 374498, 374510, 374511, 374518, 374519, 374520, 374521, 374529, 374530, 374532, 374533, 375172, 375173, 375174, 375175, 375176, 375178, 375179, 375180, 375182, 375183, 375216, 375217, 375218, 375220, 375221, 375223, 375226, 375228, 375230, 375231, 375232, 375234, 375238, 375239, 375243, 375244, 375828, 375829, 375830, 375832, 375833, 375835, 375838, 375840, 375842, 375843, 375844, 375846, 375850, 375851, 375855, 375856, 379548, 379549, 379550, 379551, 379552, 379553, 379557, 379558, 379559, 379560, 379561, 379562, 379563, 379564, 379565, 379569, 379570, 379571, 379572, 379573, 379574, 379600, 379601, 379602, 379604, 379605,
    379643, 379644, 379645, 379648, 379649, 379650, 379651, 379725, 379726, 379769
  \end{tabular}
  \caption{%
    Identified $2$-bridge alternating prime knots $16^A_i$, part 2/2.
  }%
  \label{table:Thetwobridgeknotsof16crossings2}
\end{table}


\subsection*{Acknowledgements}

I wish to thank the following persons for their assistance. Atsushi Ishii of
Osaka University composed the diagrams in the \LaTeX\/ \texttt{picture}
environment, and indeed introduced me to the $2$-bridge knots in the first
place. Ian Agol of the University of Illinois at Chicago later provided timely
advice on the form of a topological algorithm for the identification of
$2$-bridge knots. Morwen Thistlethwaite of The University of Tennessee advised
me on the modification of \textsc{Knotscape} so that it would process some of
the more knotty knots in my census, and Jon Links of The University of
Queensland provided helpful discussions.


\bibliographystyle{plain}
\bibliography{Two-Bridge_Knots}

\begin{thebibliography}{10}

\bibitem{Adams:1994}
Colin~C Adams.
\newblock {\em {The Knot Book: an Elementary Introduction to the Mathematical
  Theory of Knots}}.
\newblock Freeman, 1994.

\bibitem{DeWit:2000}
David {De Wit}.
\newblock Automatic evaluation of the {L}inks--{G}ould invariant for all prime
  knots of up to $10$ crossings.
\newblock {\em Journal of Knot Theory and its Ramifications}, 9(3):311--339,
  May 2000.
\newblock RIMS-1235, \texttt{math/9906059}.

\bibitem{DeWitKauffmanLinks:1999a}
David {De Wit}, Louis~H Kauffman, and Jon~R Links.
\newblock On the {L}inks--{G}ould invariant of links.
\newblock {\em Journal of Knot Theory and its Ramifications}, 8(2):165--199,
  March 1999.
\newblock \texttt{math/9811128}.

\bibitem{DeWitLinks:WhereLG21Fails}
David {De Wit} and Jon~R Links.
\newblock Where {$LG^{2,1}$} first fails to distinguish nonmutant prime knots.
\newblock In preparation, 2004.

\bibitem{ErnstSumners:1987}
Claus Ernst and Dewitt Sumners.
\newblock The growth of the number of prime knots.
\newblock {\em Mathematical Proceedings of the Cambridge Philosophical
  Society}, 102:303--315, 1987.

\bibitem{HosteThistlethwaiteWeeks:1998}
Jim Hoste, Morwen Thistlethwaite, and Jeff Weeks.
\newblock The first $1,701,936$ knots.
\newblock {\em The Mathematical Intelligencer}, 20(4):33--48, 1998.

\bibitem{ImafujiOchiai:2002}
Noriko Imafuji and Mitsuyuki Ochiai.
\newblock Computer aided knot theory using {M}athematica and {M}athlink.
\newblock {\em Journal of Knot Theory and its Ramifications}, 11(6):945--954,
  September 2002.

\bibitem{Ishii:2004b}
Atsushi Ishii.
\newblock Algebraic links and skein relations of the {L}inks--{G}ould
  invariant.
\newblock To appear in \emph{Proceedings of the American Mathematical Society},
  2004.

\bibitem{Kawauchi:1996}
Akio Kawauchi.
\newblock {\em A Survey of Knot Theory}.
\newblock Birkh\"{a}user Verlag, 1996.

\bibitem{Murasugi:1996}
Kunio Murasugi.
\newblock {\em Knot Theory and its Applications}.
\newblock Birkh\"{a}user, 1996.
\newblock Translated by Bohdan Kurpita.

\end{thebibliography}

\end{document}